# INTEGER ALGORITHMS TO SOLVE DIOPHANTINE LINEAR EQUATIONS AND SYSTEMS

Florentin Smarandache, Ph D
Associate Professor
Chair of Department of Math & Sciences
University of New Mexico
200 College Road
Gallup, NM 87301, USA
E-mail: smarand@unm.edu

**Abstract**: Two algorithms for solving Diophantine linear equations and five algorithms for solving Diophantine linear systems, together with properties of general and particular integer solutions, and many examples are presented in this paper.

**Contents**:

**Introduction**:

The present work includes some of the author's original researches on the integer solutions of equations and linear systems:

1. The notion of "general integer solution" of a linear equation with two unknowns is extended to linear equations with $n$ unknowns and then, to linear systems.
2. The properties of the general integer solution are determined (both of a linear equation and of a linear system).
3. Seven original integer algorithms (two for linear equations and five for linear systems) are presented. The algorithms are carefully demonstrated and an example for each of them is presented. These algorithms can be easily introduced into computer.

# INTEGER SOLUTIONS OF LINEAR EQUATIONS

Definitions and properties of the integer solutions of linear equations.

Consider the following linear equation:

(1) $$\sum_{i=1}^{n} a_i x_i = b,$$

with all $a_i \neq 0$ and $b$ in $\mathbb{Z}$.

Again, let $h \in \mathbb{N}$, and $f_i : \mathbb{Z}^h \to \mathbb{Z}$, $i = \overline{1,n}$. ($\overline{1,n}$ means: all integers from 1 to $n$).

**Definition 1.**

$x_i = x_i^0$, $i = \overline{1,n}$, is a particular integer solution of equation (1), if all $x_i^0 \in \mathbb{Z}$ and $\sum_{i=1}^{n} a_i x_i^0 = b$.

**Definition 2.**

$x_i = f_i(k_1, ..., k_h)$, $i = \overline{1,n}$, is the general integer solution of equation (1) if:

a) $\sum_{i=1}^{n} a_i f_i(k_1, ..., k_h) = b; \ \ \forall (k_1, ..., k_h) \in \mathbb{Z}^h$,

b) For any particular integer solution of equation (1), $x_i = x_i^0$, $i = \overline{1,n}$, there exist $(k_1^0, ..., k_h^0) \in \mathbb{Z}^h$ such that $x_i^0 = f_i(k_1^0, ..., k_h^0)$ for all $i = \overline{1,n}$ {i. e. any particular integer solution can be extracted from the general integer solution by parameterization}.

We will further see that the general integer solution can be expressed by linear functions.

For $1 \leq i \leq n$ we consider the functions $f_i = \sum_{j=1}^{h} c_{ij} k_j + d_i$ with all $c_{ij}, \ d_i \in \mathbb{Z}$.

**Definition 3.**

$A = (c_{ij})_{i,j}$ is the matrix associated with the general solution of equation (1).

**Definition 4.**

The integers $k_1, ..., k_s$, $1 \leq s \leq h$ are independent if all the corresponding column vectors of matrix $A$ are linearly independent.

**Definition 5.**

An integer solution is $s$-times undetermined if the maximal number of independent parameters is $s$.

**Theorem 1.** The general integer solution of equation (1) is $(n-1)$-times undetermined.

*Proof:*

We suppose that the particular integer solution is of the form:

(2) $x_i = \sum_{e=1}^{r} u_{ie} P_e + v_i, \; i = \overline{1,n}$, with all $u_{ie}, \; v_i \in \mathbb{Z}$,

$P_e$ are parameters of $\mathbb{Z}$, while $a \le r < n-1$.

Let $(x_1^0,...,x_n^0)$ be a general integer solution of equation (1) (we are not interested in the case when the equation does not have an integer solution). The solution:

$$\begin{cases} x_j = a_n k_j + x_j^0, & j = \overline{1,n-1} \\ x_n = -\left( \sum_{j=1}^{n-1} a_j k_j - x_n^0 \right) \end{cases}$$

is undetermined $(n-1)$-times (it can be easily checked that the order of the associated matrix is $n-1$). Hence, there are $n-1$ undetermined solutions. Let's consider, in the general case, a solution be undetermined $(n-1)$-times:

$x_i = \sum_{j=1}^{n-1} c_{ij} k_j + d_i, \; i = \overline{1,n}$ with all $c_{ij}, \; d_i \in \mathbb{Z}$.

Consider the case when $b = 0$.

Then

$$\sum_{i=1}^{n} a_i x_i = 0.$$

It follows:

$$\sum_{i=1}^{n} a_i x_i = \sum_{i=1}^{n} a_i \left( \sum_{j=1}^{n-1} c_{ij} k_j + d_i \right) = \sum_{i=1}^{n} a_i \sum_{j=1}^{n-1} c_{ij} k_j + \sum_{i=1}^{n} a_i d_i = 0.$$

For $k_j = 0, \; j = \overline{1,n-1}$ it follows that $\sum_{i=1}^{n} a_i d_i = 0$.

For $k_{j_0} = 1$ and $k_j = 0, \; j \ne j_0$, it follows that $\sum_{i=1}^{n} a_i c_{ij_0} = 0$.

Let's consider the homogenous linear system of $n$ equations with $n$ unknowns:

$$\begin{cases} \sum_{i=1}^{n} x_i c_{ij} = 0, & j = \overline{1,n-1} \\ \sum_{i=1}^{n} x_i d_i = 0 \end{cases}$$

which, obviously, has the solution $x_i = a_i, \; i = \overline{1,n}$ different from the trivial one. Hence the determinant of the system is zero, i.e., the vectors $c_j = (c_{1j},...,c_{nj})^t, \; j = \overline{1,n-1}$, $D = (d_1,...,d_n)^t$ are linearly dependent.

But the solution being $(n-1)$-times undetermined it shows that $c_j,\ j=\overline{1,n-1}$ are linearly independent. Then $(c_1,...,c_{n-1})$ determines a free sub-module $\mathbb{Z}$ of order $n-1$ in $\mathbb{Z}_n$ of solutions for the given equation.

Let's see what can we obtain from (2). We have:

$$0=\sum_{i=1}^{n}a_ix_i=\sum_{i=1}^{n}a_i\left(\sum_{e=1}^{r}u_{ie}P_e+v_i\right).$$

As above, we obtain:

$$\sum_{i=1}^{n}a_iv_i=0 \text{ and } \sum_{e=1}^{r}a_iu_{ie_0}=0$$

similarly, the vectors $U_h=(u_{1h},...,u_{nh})$ are linearly independent, $h=\overline{1,r}$, $U_h,\ h=\overline{1,r}$ are $V=(v_1,...,v_n)$ particular integer solutions of the homogenous linear equation.

**Sub-case (a1)**

$U,h=\overline{1,r}$ are linearly dependent. This gives $\{U_1,...,U_r\}=$ the free sub-module of order $r$ in $\mathbb{Z}^n$ of solutions of the equation. Hence, there are solutions from $\{V_1,...,V_{n-1}\}$ which are not from $\{U_1,...,U_r\}$; this contradicts the fact that (2) is the general integer solution.

**Sub-case (a2)**

$U_h,\ h=\overline{1,r},\ V$ are linearly independent. Then $\{U_1,...,U_r\}+V$ is a linear variety of the dimension $<n-1=\dim\{V_1,...,V_{n-1}\}$ and the conclusion can be similarly drawn.

Consider the case when $b\neq 0$. So, $\sum_{i=1}^{n}a_ix_i=b$.

Then:

$$\sum_{i=1}^{n}a_i\left(\sum_{j=1}^{n-1}c_{ij}k_j+d_i\right)=\sum_{j=1}^{n-1}\left(\sum_{i=1}^{n}a_ic_{ij}\right)k_j+\sum_{i=1}^{n}a_id_i=b;\ \ \forall(k_1,...,k_{n-1})\in\mathbb{Z}^{n-1}.$$

As in the previous case, we obtain $\sum_{i=1}^{n}a_id_i=b$ and $\sum_{i=1}^{n}a_ic_{ij}=0, \quad \forall\, j=\overline{1,n-1}$.

The vectors $c_j=(c_{ij},...,c_{nj})^t,\ j=\overline{1,n-1}$, are linearly independent because the solution is undetermined $(n-1)$-times.

Conversely, if $c_1,...,c_{n-1},\ D$ (where $D=(d_1,...,d_n)^t$) were linearly dependent, it would mean that $D=\sum_{j=1}^{n-1}s_jc_j$ with all $s_j$ scalar; it would also mean that

$$b=\sum_{i=1}^{n}a_id_i=\sum_{i=1}^{n}a_i\left(\sum_{j=1}^{n-1}s_jc_{ij}\right)=\sum_{j=1}^{n-1}s_j\left(\sum_{i=1}^{n}a_ic_{ij}\right)=0.$$

This is impossible.

(3) Then $\{c_1,...,c_{n-1}\}+D$ is a linear variety.

Let us see what we can obtain from (2). We have:

$$b = \sum_{i=1}^{n} a_i x_i = \sum_{i=1}^{n} a_i \left( \sum_{e=1}^{r} u_{ie} P_e + v_i \right) = \sum_{e=1}^{r} \left( \sum_{i=1}^{n} a_i u_{ie} \right) P_e + \sum_{i=1}^{n} a_i v_i$$

and, similarly: $\sum_{i=1}^{n} a_i v_i = b$ and $\sum_{i=1}^{n} a_i u_{ie} = 0, \quad \forall\, e = \overline{1,r}$, respectively. The vectors $U_e = (u_{1e},...,u_{ne})^t$, $e = \overline{1,r}$ are linearly independent because the solution is undetermined $r$-times.

A procedure like that applied in (3) shows that $U_1,...,U_r, V$ are linearly independent, where $V = (v_1,...,v_n)^t$. Then $\{U_1,...,U_r\} + V$ = a linear variety = free sub-module of order $r < n-1$. That is, we can find vectors from $\{c_1,...,c_{n-1}\} + D$ which are not from $\{U_1,...,U_r\} + V$, contradicting the "general" characteristic of the integer number solution. Hence, the general integer solution is undetermined $(n-1)$-times.

**Theorem 2.** The general integer solution of the homogeneous linear equation $\sum_{i=1}^{n} a_i x_i = 0$ (all $a_i \in \mathbb{Z} \setminus \{0\}$) can be written under the form:

$$(4) \qquad x_i = \sum_{j=1}^{n-1} c_{ij} k_j, \;\; i = \overline{1,n}$$

(with $d_1 = ... = d_n = 0$).

**Definition 6.** This is called the standard form of the general integer solution of a homogeneous linear equation.

*Proof:*

We consider the general integer solution under the form:

$$x_i = \sum_{j=1}^{n-1} c_{ij} P_j + d_i, \;\; i = \overline{1,n}$$

with not all $d_i = 0$. We'll show that it can be written under the form (4). The homogeneous equation has the trivial solution $x_i = 0,\; i = \overline{1,n}$. There is $(p_1^0,...,p_{n-1}^0) \in \mathbb{Z}^{n-1}$ such that $\sum_{j=1}^{n-1} c_{ij} p_j^0 + d_i = 0, \;\; \forall i = \overline{1,n}$.

Substituting: $P_j = k_j + p_j, \;\; j = \overline{1,n-1}$ in the form shown at the beginning of the demonstration, we will obtain form (4). We have to mention that the substitution does not diminish the degree of generality as $P_j \in \mathbb{Z} \Leftrightarrow k_j \in \mathbb{Z}$ because $j = \overline{1,n-1}$.

**Theorem 3.** The general integer solution of a non-homogeneous linear equation is equal to the general integer solution of its associated homogeneous linear equation plus any particular integer solution of the non-homogeneous linear equation.

*Proof:*

Let's consider that $x_i = \sum_{j=1}^{n-1} c_{ij} k_j, \quad i = \overline{1,n}$, is the general integer solution of the associated homogeneous linear equation and, again, let $x_i = v_i, \ i = \overline{1,n}$, be a particular integer solution of the non-homogeneous linear equation. Then $x_i = \sum_{j=1}^{n-1} c_{ij} k_j + v_i, \quad i = \overline{1,n}$, is the general integer solution of the non-homogeneous linear equation.

Actually, $\sum_{i=1}^{n} a_i x_i = \sum_{i=1}^{n} a_i \left( \sum_{j=1}^{n-1} c_{ij} k_j + v_i \right) = \sum_{i=1}^{n} a_i \left( \sum_{j=1}^{n-1} c_{ij} k_j \right) + \sum_{i=1}^{n} a_i v_i = b$;

if $x_i = x_i^0, \ i = \overline{1,n}$, is a particular integer solution of the non-homogeneous linear equation, then $x_i = x_i - v_i, \ i = \overline{1,n}$, is a particular integer solution of the homogeneous linear equation: hence, there is $(k_1^0, ..., k_{n-1}^0) \in \mathbb{Z}^{n-1}$ such that

$$\sum_{j=1}^{n-1} c_{ij} k_j^0 = x_i^0 - v_i, \ \forall \ i = \overline{1,n},$$

i.e.:

$$\sum_{j=1}^{n-1} c_{ij} k_j^0 + v_i = x_i^0, \ \forall \ i = \overline{1,n},$$

which was to be proven.

**Theorem 4.** If $x_i = \sum_{j=1}^{n-1} c_{ij} k_j, \ i = \overline{1,n}$ is the general integer solution of a homogeneous linear equation $\left(c_{ij}, ..., c_{nj}\right) \sim 1 \ \forall j = \overline{1, n-1}$.

The demonstration is done by reduction ad absurdum. If $\exists j_0, \ 1 \le j_0 \le n-1$ such that $\left(c_{ij_0}, ..., c_{nj_0}\right) \sim d_{j_0} \ne \pm 1$, then $c_{ij_0} = c'_{ij_0} d_{ij_0}$ with $\left(c'_{ij_0}, ..., c'_{nj_0}\right) \sim 1, \quad \forall \ i = \overline{1,n}$.

But $x_i = c'_{ij_0}, \quad i = \overline{1,n}$, represents a particular integer solution as

$$\sum_{i=1}^{n} a_i x_i = \sum_{i=1}^{n} a_i c'_{ij_0} = 1 / d_{j_0} \cdot \sum_{i=1}^{n} a_i c_{ij_0} = 0$$

(because $x_i = c_{ij_0}, \ i = \overline{1,n}$ is a particular integer solution from the general integer solution by introducing $k_{j_0} = 1$ and $k_j = 0, \ j \ne j_0$. But the particular integer solution $x_i = c'_{ij_0}, \ i = \overline{1,n}$, cannot be obtained by introducing integer number parameters (as it should) from the general integer solution, as from the linear system of $n$ equations and $n-1$ unknowns, which is compatible. We obtain:

$$x_i = \sum_{\substack{j=1 \\ j \ne j_0}}^{n} c_{ij} k_j + c'_{ij_0} d_{j_0} k_{j_0} = c'_{ij_0}, \quad i = \overline{1,n}.$$

Leaving aside the last equation – which is a linear combination of other $n-1$ equations – a Kramerian system is obtained, as follows:

$$k_{j_0} = \frac{\begin{vmatrix} c_{11}......c'_{ij_0}.......c_{1,n-1} \\ \vdots \\ c_{n-1,1}...c'_{n-1j_0}...c_{n-1n-1} \end{vmatrix}}{\begin{vmatrix} c_{11}.......c'_{ij_0}d_{j_0}......c_{1,n-1} \\ \vdots \\ c_{n-1,1}...c'_{n-1j_0}d_{j_0}...c_{n-1n-1} \end{vmatrix}} = \frac{1}{d_{j_0}} \notin \mathbb{Z}$$

Therefore the assumption is false (end of demonstration).

**Theorem 5.** Considering the equation (1) with $\left(a_1,...,a_n\right) \sim 1,\ \ b = 0$ and the general integer solution $x_i = \sum_{j=1}^{n-1} c_{ij}k_j,\ \ i = \overline{1,n}$, then

$$\left(a_1,...,a_{i-1},a_{i+1},...,a_n\right) \sim \left(c_{i1},...,c_{in-1}\right),\ \ \forall i = \overline{1,n}.$$

*Proof:*

The demonstration is done by double divisibility.

Let's consider $i_0,\ 1 \le i_0 \le n$ arbitrary but fixed. $x_{i_0} = \sum_{j=1}^{n-1} c_{i_0 j}k_j$. Consider the equation $\sum_{\mathrm{i}\ne \mathrm{i}_0} a_i x_i = -a_{i_0}x_{i_0}$. We have shown that $x_i = c_{ij},\ \ i = \overline{1,n}$ is a particular integer solution irrespective of $j,\ \ a \le j \le n-1$.

The equation $\sum_{\mathrm{i}\ne \mathrm{i}_0} a_i x_i = -a_{i_0}c_{i_0 j}$ obviously, has the integer solution $x_i = c_{ij},\ \ i \ne i_0$. Then $\left(a_1,...,a_{i_0-1},a_{i_0+1},...,a_n\right)$ divides $-a_{i_0}c_{i_0 j}$ as we have assumed, it follows that $\left(a_1,...,a_n\right) \sim 1$, and it follows that $\left(a_1,...,a_{i_0-1},a_{i_0+1},...,a_n\right) \mid c_{i_0 j}$ irrespective of $j$. Hence $\left(a_1,...,a_{i_0-1},a_{i_0+1},...,a_n\right) \mid \left(c_{i_0 1},...,c_{i_0 n-1}\right),\ \ \forall i = \overline{1,n}$, and the divisibility in one sense was proven.

Inverse divisibility:

Let us suppose the contrary and consider that $\exists i_1 \in \overline{1,n}$ for which $\left(a_1,...,a_{i_1-1},a_{i_1+1},...,a_n\right) \sim d_{i_1 1} \ne d_{i_1 2} \sim \left(c_{i_1 1},...,c_{i_1 n-1}\right)$; we have considered $d_{i_1 1}$ and $d_{i_1 2}$ without restricting the generality. $d_{i_1 1} \mid d_{i_1 2}$ according to the first part of the demonstration. Hence, $\exists d \in \mathbb{Z}$ such that $d_{i_1 2} = d \cdot d_{i_1 1},\ \ |d| \ne 1$.

$$x_{i_1} = \sum_{j=1}^{n-1} c_{i_1 j}k_j = d \cdot d_{i_1 1}\sum_{j=1}^{n-1} c'_{i_1 j}k_j\ ;$$

$$\sum_{i=1}^{n} a_i x_i = 0 \Rightarrow \sum_{i \ne i_1}^{n} a_i x_i = -a_{i_1}x_{i_1} \sum_{i \ne i_1} a_i x_i = -a_{i_1} d \cdot d_{i_1 1}\sum_{j=1}^{n-1} c'_{i_1 j}k_j\ ,$$

where $\left(c_{i_1 1},...,c_{i_1 n-1}\right) \sim 1$.

The non-homogeneous linear equation $\sum_{i \neq i_1} a_i x_i = -a_{i_1} d_{i_1 1}$ has the integer solution because $a_{i_1} d_{i_1 1}$ is divisible by $\left(a_1,...,a_{i_1-1},a_{i_1+1},...,a_n\right)$. Let's consider that $x_i = x_i^0,\ i \neq i_1$, is its particular integer solution. It follows that the equation $\sum_{i=1}^{n} a_i x_i = 0$ the particular solution $x_i = x_i^0,\ i \neq i_1,\ \ x_{i_1} = d_{i_1}$, which is written as (5). We'll show that (5) cannot be obtained from the general solution by integer number parameters:

$$(6) \qquad \begin{cases} \sum_{j=1}^{n-1} c_{ij} k_j = x_i^0, \quad i \neq i_1 \\ d \cdot d_{i_1 1} \sum_{j=1}^{n-1} c_{ij} k_j = d_{i_1 1} \end{cases}$$

But the equation (6) does not have an integer solution because $d \cdot d_{i_1 1} \nmid d_{i_1 1}$ thus, contradicting, the "general" characteristic of the integer solution.

As a conclusion we can write:

**Theorem 6**. Let's consider the homogeneous linear equation $\sum_{i=1}^{n} a_i x_i = 0$, with all $a_i \in \mathbb{Z} \setminus \{0\}$ and $\left(a_1,...,a_n\right) \sim 1$.

Let $x_i = \sum_{j=1}^{h} c_{ij} k_j,\ i = \overline{1,n}$, with all $c_{ij} \in \mathbb{Z}$, all $k_j$ integer parameters and let's consider $h \in \mathbb{N}$ be a general integer solution of the equation. Then,

1) the solution is undetermined $(n-1)$-times;
2) $\forall\, j = \overline{1,n-1}$ we have $\left(c_{1j},...,c_{nj}\right) \sim 1$;
3) $\forall\, i = \overline{1,n}$ we have $\left(c_{i1},...,c_{in-1}\right) \sim \left(a_1,...,a_{i-1},a_{i+1},...,a_n\right)$.

The proof results from theorems 1,4 and 5.

**Note 1**. The only equation of the form (1) that is undetermined $n$-times is the trivial equation $0 \cdot x_1 + ... + 0 \cdot x_n = 0$.

**Note 2**. The converse of theorem 6 is not true.

Counterexample:

$$(7) \qquad \begin{cases} x_1 = -k_1 + k_2 \\ x_2 = 5k_1 + 3k_2 \\ x_3 = 7k_1 - k_2; \quad k_1, k_2 \in \mathbb{Z} \end{cases}$$

is not the general integer solution of the equation

$$(8) \qquad -13x_1 + 3x_2 - 4x_3 = 0$$

although the solution (7) verifies the points 1), 2) and 3) of theorem 6. (1, 7, 2) is the particular integer solution of (8) but cannot be obtained by introducing integer number parameters in (7) because from

$$\begin{cases} -k_1 + k_2 = 1 \\ 5k_1 + 3k_2 = 7 \\ 7k_1 - k_2 = 2 \end{cases}$$

it follows that $k = \frac{1}{2} \notin \mathbb{Z}$ and $k = \frac{3}{2} \notin \mathbb{Z}$ (unique roots).

**REFERENCE**

[1] Smarandache, Florentin – Whole number solution of linear equations and systems – diploma thesis work, 1979, University of Craiova (under the supervision of Assoc. Prof. Dr. Alexandru Dincă)

# AN INTEGER NUMBER ALGORITHM TO SOLVE LINEAR EQUATIONS

An algorithm is given that ascertains whether a linear equation has integer number solutions or not; if it does, the general integer solution is determined.

**Input**

A linear equation $a_1x_1 + ... + a_nx_n = b$, with $a_i, b \in \mathbb{Z}$, $x_i$ being integer number unknowns, $i = \overline{1,n}$, and not all $a_i = 0$.

**Output**

Decision on the integer solution of this equation; and if the equation has solutions in $\mathbb{Z}$, its general solution is obtained.

**Method**

*Step 1*. Calculate $d = (a_1, ..., a_n)$.

*Step 2*. If $d / b$ then "the equation has integer solution"; go on to Step 3. If $d \nmid b$ then "the equation does not have integer solution"; stop.

*Step 3*. Consider $h := 1$. If $|d| \neq 1$, divide the equation by $d$; consider $a_i := a_i / d,\ i = \overline{1,n},\ b := b / d$.

*Step 4*. Calculate $a = \min\limits_{a_s \neq 0} |a_s|$ and determine an $i$ such that $a_i = a$.

*Step 5*. If $a \neq 1$ then go to Step 7.

*Step 6*. If $a = 1$, then:

(A) $x_i = -(a_1x_1 + ... + a_{i-1}x_{i-1} + a_{i+1}x_{i+1} + ... + a_nx_n - b) \cdot a_i$

(B) Substitute the value of $x_i$ in the values of the other determined unknowns.

(C) Substitute integer number parameters for all the variables of the unknown values in the right term: $k_1, k_2, ..., k_{n-2}$, and $k_{n-1}$ respectively.

(D) Write, for your records, the general solution thus determined; stop.

*Step7*. Write down all $a_j,\ j \neq i$ and under the form:

$$a_j = a_iq_j + r_j$$

$$b = a_iq + r \text{ where } q_j = \left[\frac{a_j}{a_i}\right],\ q = \left[\frac{b}{a_i}\right].$$

*Step 8*. Write $x_i = -q_1x_1 - ... - q_{i-1}x_{i-1} - q_{i+1}x_{i+1} - ... - q_nx_n + q - t_h$. Substitute the value of $x_i$ in the values of the other determined unknowns.

*Step 9*. Consider

$$\begin{cases} a_1 := r_1 \\ : \\ a_{i-1} := r_{i-1} \\ a_{i+1} := r_{i+1} \\ : \\ a_n := r_n \end{cases} \quad \text{and} \quad \begin{cases} a_i := -a_i \\ b := r \\ x_i := t_h \\ h := h+1 \end{cases}$$

and go back to Step 4.

**Lemma 1.** The previous algorithm is finite.

*Proof:*

Let's $a_1x_1 + ... + a_nx_n = b$ be the initial linear equation, with not all $a_i = 0$; check for $\min_{a_s \neq 0} |a_s| = a_1 \neq 1$ (if not, it is renumbered). Following the algorithm, once we pass from this initial equation to a new equation: $a_1^{'}x_1 + a_2^{'}x_2 + ... + a_n^{'}x_n = b'$, with $|a_1^{'}| < |a_i|$ for $i = \overline{2,n}$, $|b'| < |b|$ and $a_1^{'} = -a_1$.

It follows that $\min_{a_s^{'} \neq 0} |a_s^{'}| < \min_{a_s \neq 1} |a_s|$. We continue similarly and after a finite number of steps we obtain, at Step 4, $a := 1$ (the actual $a$ is always smaller than the previous $a$, according to the previous note) and in this case the algorithm terminates.

**Lemma 2.** Let the linear equation be:

(25) $a_1x_1 + a_2x_2 + ... + a_nx_n = b$, with $\min_{a_s \neq 0} |a_s| = a_1$ and the equation

(26) $-a_1t_1 + r_2x_2 + ... + r_nx_n = r$, with $t_1 = -x_1 - q_2x_2 - ... - q_nx_n + q$, where $r_i = a_i - a_iq_i$, $i = \overline{2,n}$, $r = b - a_1q$ while $q_i = \left[\frac{a_i}{a}\right]$, $r = \left[\frac{b}{a_1}\right]$. Then $x_1 = x_1^0$, $x_2 = x_2^0, \ldots, x_n = x_n^0$ is a particular solution of equation (25) if and only if $t_1 = t_1^0 = -x_1 - q_2x_2^0 - ... - q_nx_n^0 + q$, $x_2, ..., x_n = x_n^0$ is a particular solution of equation (26).

*Proof:*

$x_1 = x_1^0$, $x_2 = x_2^0, \ldots, x_n = x_n^0$, is a particular solution of equation (25)$\Leftrightarrow$
$a_1x_1^0 + a_2x_2^0 + ... + a_nx_n^0 = b \Leftrightarrow a_1x_1^0 + (r_2 + a_1q_2)x_2^0 + ... + (r_n + a_1q_n)x_n^0 = a_1q + r \Leftrightarrow$
$r_2x_2^0 + ... + r_nx_n^0 - a_1(-x_1^0 - q_2x_2^0 - ... - q_nx_n^0 + q) = r \Leftrightarrow -a_1t_1^0 + r_2x_2^0 + ... + r_nx_n^0 = r \Leftrightarrow$
$\Leftrightarrow t_1 = t_1^0, x_2 = x_2^0, ..., x_n = x_n^0$ is a particular solution of equation (26).

**Lemma 3.** $x_i = c_{i1}k_1 + ... + c_{in-1}k_{n-1} + d_i$, $i = \overline{1,n}$, is the general solution of equation (25) if and only if

(28) $t_1 = -(c_{11} + q_2c_{21} + ... + q_nc_{n1})k_1 - ... - (c_{1n-1} + q_2c_{2n-1} + ... + q_nc_{nn-1})k_n -$
$-(d_1 + q_2d_2 + ... + q_nd_n) + q$,
$x_j = c_{1j1}k_1 + ... + c_{jn-1}k_{n-1} + d_j$, $j = \overline{2,n}$

is a general solution for equation (26).

*Proof:*

$t_1 = t_1^0 = -x_1^0 - q_2 x_2^0 - \ldots - q_n x_n^0 + q,\ x_2 = x_2^0, \ldots, x_n = x_n^0$ is a particular solution of the equation (25) $\Leftrightarrow$ $x_1 = x_1^0,\ x_2 = x_2^0, \ldots, x_n = x_n^0$ is a particular solution of equation (26) $\Leftrightarrow$ $\exists k_1 = k_1^0 \in \mathbb{Z}, \ldots, k_n = k_n^0 \in \mathbb{Z}$ such that

$$x_i = c_{i1}k_1^0 + \ldots + c_{in-1}k_{n-1}^0 + d_i = x_i^0,\ i = \overline{1,n} \Leftrightarrow \exists k_1 = k_1^0 \in \mathbb{Z}, \ldots, k_n = k_n^0 \in \mathbb{Z},$$

such that

$$x_i = c_{i1}k_1^0 + \ldots + c_{in-1}k_{n-1}^0 + d_i = x_i^0,\ \ i = \overline{2,n},$$

and

$$t_1 = -(c_{11} + q_2 c_{21} + \ldots + q_n c_{n1})k_1^0 - \ldots - (c_{1n-1} + q_2 c_{2n-1} + \ldots + q_n c_{nn-1})k_{n-1}^0 -$$
$$(d_1 + q_2 d_2 + \ldots + q_n d_n) + q = -x_1^0 - q_2 x_2^0 - \ldots - q_n x_n^0 + q = t_1^0$$

**Lemma 4.** The linear equation

(29) $\quad a_1 x_1 + a_2 x_2 + \ldots + a_n x_n = b$ with $|a_1| = 1$ has the general solution:

(30) $$\begin{cases} x_1 = -(a_2 k_2 + \ldots + a_n k_n - b)a_1 \\ x_i = k_i \in \mathbb{Z} \\ i = \overline{2,n} \end{cases}$$

*Proof:*

Let's consider $x_1 = x_1^0,\ x_2 = x_2^0, \ldots, x_n = x_n^0$, a particular solution of equation (29). $\exists k_2 = x_2^0,\ k_n = x_n^0$, such that $x_1 = \left(-a_2 x_2^0 + \ldots + a_n x_n^0 - b\right)a_1 = x_1^0,\ \ x_2 = x_2^0, \ldots, x_n = x_n^0$.

**Lemma 5.** Let's consider the linear equation $a_1 x_1 + a_2 x_2 + \ldots + a_n x_n = b$, with $\min\limits_{a_s \neq 0} |a_s| = a_1$ and $a_i = a_1 q_i,\ \ i = \overline{2,n}$.

Then, the general solution of the equation is:

$$\begin{cases} x_1 = -(q_2 k_2 + \ldots + q_n k_n - q) \\ x_i = k_i \in \mathbb{Z} \\ i = \overline{2,n} \end{cases}$$

*Proof:*

Dividing the equation by $a_1$ the conditions of Lemma 4 are met.

**Theorem of Correctness.** The preceding algorithm calculates correctly the general solution of the linear equation $a_1 x_1 + \ldots + a_n x_n = b$, with not all $a_i = 0$.

*Proof:*

The algorithm is finite according to Lemma 1. The correctness of steps 1, 2, and 3 is obvious. At step 4 there is always $\min\limits_{a_s \neq 0} |a_s|$ as not all $a_i = 0$. The correctness of sub-step 6 A) results from Lemmas 4 and 5, respectively. This algorithm represents a method of obtaining the general solution of the initial equation by means of the general solutions of the linear equation obtained after the algorithm was followed several times (according to Lemmas 2 and 3); from Lemma 3, it follows that to obtain the general solution of the initial linear equation is equivalent to calculate the general solution of an equation at step 6 A), equation whose general solution is given in algorithm (according to Lemmas 4 and 5). The Theorem of correctness has been fully proven.

**Note**. At step 4 of the algorithm we consider $a := \min_{a_s \neq 0} |a_s|$ such that the number of iterations is as small as possible. The algorithm works if we consider $a := |a_i| \neq \max_{s=\overline{1,n}} |a_s|$ but it takes longer. The algorithm can be introduced into a computer program.

**Application**

Calculate the integer solution of the equation:

$$6x_1 - 12x_2 - 8x_3 + 22x_4 = 14 .$$

*Solution*

The previous algorithm is applied.

1. $(6,-12,-8,22) = 2$
2. $2 \mid 14$ therefore the solution of the equation is in $\mathbb{Z}$ .
3. $h := 1$; $|2| \neq 1$; dividing the equation by 2 we obtain:

$3x_1 = 6x_2 - 4x_3 + 11x_4 = 7$ .

4. $a := \min\{|3|,\ |-6|,\ |-4|,\ |11|\} = 3,\ i = 1$
5. $a \neq 1$
7. $-6 = 3 \cdot (-2) + 0$

$-4 = 3 \cdot (-2) + 2$

$11 = 3 \cdot 3 + 2$

$7 = 3 \cdot 2 + 1$

8. $x_1 = 2x_2 + 2x_3 - 3x_4 + 2 - t_1$
9.

$a_2 := 0 \qquad a_1 := -3$

$a_3 := 2 \qquad b := 1$

$a_4 := 2 \qquad x_1 := t_1$

$h := 2$

4. We have a new equation:

$-3t_1 - 0 \cdot x_2 + 2x_3 + 2x_4 = 1$

$a := \min\{|-3|,|2|,|2|\}$ and

$i = 3$

5. $a \neq 1$
7. $-3 = 2 \cdot (-2) + 1$

$0 = 2 \cdot 0 + 0$

$2 = 2 \cdot 1 + 0$

$1 = 2 \cdot 0 + 0$

8. $x_3 = 2t_1 + 0 \cdot x_2 - x_4 + 0 - t_2$. Substituting the value of $x_3$ in the value determined for $x_1$ we obtain: $x_1 = 2x_2 - 5x_4 + 3t_1 - 2t_2 + 2$

9. $a_1 := 1$ $a_3 := -2$
$a_2 := 0$ $b := 1$
$a_4 := 0$ $x_3 := t_2$
$h := 3$

4. We have obtained the equation:
$1 \cdot t_2 + 0 \cdot x_2 - 2 \cdot t_2 + 0 \cdot x_4 = 1$, $a = 1$, and $i = 1$

6. (A) $t_1 = -(0 \cdot x_2 - 2t_2 + 0 \cdot x_4 - 1) \cdot 1 = 2t_2 + 1$

(B) Substituting the value of $t_1$ in the values of $x_1$ and $x_3$ previously determined, we obtain:

$$x_1 = 2x_2 - 5x_4 + 4t_2 + 5 \text{ and}$$
$$x_3 = -x_4 + 3t_2 + 2$$

(C) $x_2 := k_1$, $x_4 := k_2$, $t_2 := k_3$, $k_1, k_2, k_3 \in \mathbb{Z}$

(D) The general solution of the initial equation is:

$$x_1 = 2k_1 - 5k_2 + 4k_3 + 5$$
$$x_2 = k_1$$
$$x_3 = -k_2 + 3k_3 + 2$$
$$x_4 = k_2$$

$k_1, k_2, k_3$ are parameters $\in \mathbb{Z}$

## REFERENCE

[1] Smarandache, Florentin – Whole number solution of equations and systems of equations – part of the diploma thesis, University of Craiova, 1979.

# ANOTHER INTEGER NUMBER ALGORITHM TO SOLVE LINEAR EQUATIONS (USING CONGRUENCES)

In this section is presented a new integer number algorithm for linear equation. This algorithm is more “rapid” than W. Sierpinski’s presented in [1] in the sense that it reaches the general solution after a smaller number of iterations. Its correctness will be thoroughly demonstrated.

## ANOTHER INTEGER NUMBER ALGORITHM TO SOLVE LINEAR EQUATIONS

Let’s us consider the equation (1); (the case $a_i, b \in \mathbb{Q},\ i = \overline{1,n}$ is reduced to the case (1) by reducing to the same denominator and eliminating the denominators). Let $d = (a_1, ..., a_n)$. If $d \mid b$ then the equation does not have integer solutions, while if $d \nmid b$ the equation has integer solutions (according to a well-known theorem from the number theory).

If the equation has solutions and $d \neq$ we divide the equation by $d$. Then $d = 1$ (we do not make any restriction if we consider the maximal co-divisor positive).

Also,

(a) If all $a_i$ the equation is trivial; it has the general integer solution $x_i = k_i \in \mathbb{Z},\ i = \overline{1,n}$, when $b = 0$ (the only case when the general solution is $n$-times undetermined) and does not have solution when $b \neq 0$.

(b) If $\exists i,\ 1 \le i \le n$ such that $a_i = \pm 1$ then the general integer solution is: $x_i = -a_i \left( \sum_{\substack{j=1 \\ j \neq i}}^{n} a_j k_j - b \right)$ and $x_s = k_s \in \mathbb{Z},\ s \in \{1, ..., n\} \setminus \{i\}$

The proof of this assertion was given in [4]. All these cases are trivial, therefore we will leave them aside. The following algorithm can be written:

**Input**

A linear equation:

(2) $$\sum_{i=1}^{n} a_i x_i = b,\ a_i, b \in \mathbb{Z},\ \ a_i \neq \pm 1,\ i = \overline{1,n},$$

with not all $a_i = 0$ and $(a_1, ..., a_n) = 1$.

**Output**

The integer general solution of the equation.

**Method**

1. $h := 1,\ p := 1$

2. Calculate $\min\limits_{1\le i,j\le n}\left\{|r|,\ r\equiv a_i(\mathrm{mod}\, a_j),\ |r|<|a_j|\right\}$ and determine $r$ and the pair $(i,j)$ for which this minimum can be obtained (when there are more possibilities we have to choose one of them).

3. If $|r|\neq$ go to step 4.

If $|r|=1$, then

$$\begin{cases} x_i := r\left(-a_j t_h - \sum\limits_{\substack{s=1\\ s\notin\{i,j\}}}^{n} a_s x_s + b\right) \\ x_j := r\left(a_i t_h + \dfrac{a_i - r}{a_j}\cdot \sum\limits_{\substack{s=1\\ s\notin\{i,j\}}}^{n} a_s x_s + \dfrac{r-a_i}{a_j} b\right) \end{cases}$$

(A) Substitute the values thus determined of these unknowns in all the statements $(p),\ p=1,2,...$ (if possible).

(B) From the last relation $(p)$ obtained in the algorithm substitute in all relations: $(\bar{p}-1),(\bar{p}-2),...,(1)$

(C) Every statement, starting in order from $(\bar{p}-1)$ should be applied the same procedure as in (B): then $(\bar{p}-2),...,(3)$ respectively.

(D) Write the values of the unknowns $x_i,\ i=\overline{1,n}$, from the initial equation (writing the corresponding integer number parameters from the right term of these unknowns with $k_1,...,k_{n-1}$), STOP.

4. Write the statement $(p)$: $x_j = t_h - \dfrac{a_i - r}{a_j} x_i$

5. Assign $\quad x_j := t_h \qquad h := h+1$

$\qquad\qquad\quad a_i := r \qquad p := p+1$

The other coefficients and variables remain unchanged go back to step 2.

**The Correctness of the Algorithm**

Let us consider linear equation (2). Under these conditions, the following properties exist:

**Lemma 1.** The set $M=\left\{\ |r|,\ r\equiv a_i(\mathrm{mod}\, a_j),\ 0<|r|<|a_j|\right\}$ has a minimum.

*Proof:*

Obviously $M\subset\mathbb{N}^*$ and $M$ is finite because the equation has a finite number of coefficients: $n$, and considering all the possible combinations of these, by twos, there is the maximum $AR_n^2$ (arranged with repetition) = $n^2$ elements.

Let us show, by reduction ad absurdum, that $M \neq \emptyset$.

$M \neq \emptyset \Leftrightarrow a_i \equiv 0(\bmod a_j)\ \forall i, j = \overline{1,n}$. Hence $a_j \equiv 0(\bmod a_i)\ \forall i, j = \overline{1,n}$. Or this is possible only when $|a_i| = |a_j|,\ \forall i, j = \overline{1,n}$, which is equivalent to $(a_1,..,a_n) = a_i,\ \forall i \in \overline{1,n}$. But $(a_1,..,a_n) = 1$ are a restriction from the assumption. It follows that $|a_i| = \overline{1,n},\ \forall i \in \overline{1,n}$ a fact which contradicts the other restrictions of the assumption.

$M \neq 0$ and finite, it follows that $M$ has a minimum.

**Lemma 2.** If $|r| = \min\limits_{1\le i,j\le n} M$ , then $|r| < |a_i|,\ \forall i \in \overline{1,n}$.

*Proof:*

We assume conversely, that $\exists i_0,\ 1 \le i_0 \le n$ such that $|r| \ge |a_{i_0}|$.

Then $|r| \ge \min\limits_{1\le j\le n}\left\{\ |a_j|\right\} = |a_{j_0}| \neq 1,\ 1 \le j_0 \le n$. Let $a_{p_0},\ 1 \le p_0 \le n$, such that $|a_{p_0}| > |a_{j_0}|$ and $a_{p_0}$ is not divided by $a_j^0$.

There is a coefficient in the equation, $|a_{j_0}|$ which is the minimum and the coefficients are not equal among themselves (conversely, it would mean that $(a_1,..,a_n) = a_1 = \pm 1$ which is against the hypothesis and, again, of the coefficients whose absolute value is greater that $|a_{ij_0}|$ not all can be divided by $a_{j_0}$ (conversely, it would similarly result in $(a_1,..,a_n) = a_{j_0} \neq \pm 1$.

We write $\left[a_{p_0} / a_{j_0}\right] = q_0 \in \mathbb{Z}$ (integer portion), and $r = a_{p_0} - q_0 a_{j_0} \in \mathbb{Z}$. We have $a_{p_0} \equiv r_0(\bmod a_{j_0})$ and $0 < |r_0| < |a_{j_0}| < |a_{i_0}| \le |r|$. Thus, we have found an $r_0$ which $|r_0| < |r|$ contradicts the definition of minimum given to $|r|$.

Thus $|r| < |a_i|,\ \forall i \in \overline{1,n}$.

**Lemma 3.** If $|r| = \min M = 1$ for the pair of indices $(i, j)$, then:

$$
\begin{cases}
x_i = r\left(-a_j t_h - \sum\limits_{\substack{s=1\\ s\notin\{i,j\}}}^{n} a_s k_s + b\right)\\
x_j = r\left(a_i t_h + \dfrac{a_i - r}{a_j} \cdot \sum\limits_{\substack{s=1\\ s\notin\{i,j\}}}^{n} a_s k_s + \dfrac{r - a_i}{a_j} b\right)\\
x_s = k_s \in \mathbb{Z},\ s \in \{1,...,n\} \setminus \{i, j\}
\end{cases}
$$

is the general integer solution of equation (2).

*Proof:*

Let $x_e = x_e^0,\ e=\overline{1,n}$, be a particular integer solution of equation (2). Then $\exists k_s = x_s^0 \in \mathbb{Z},\ s \in \{1,...,n\} \setminus \{i,j\}$ and $t_h = x_j^0 + \frac{a_i - r}{a_j} x_i^0 \in \mathbb{Z}$ (because $a_i - r = Ma_j$) such that:

$$x_i = r - a_j\left(x_j^0 + \frac{a_i - r}{a_j} x_i\right) - \sum_{\substack{s=1\\ s\notin\{i,j\}}}^{n} a_s x_s^0 + b = x_i^0$$

$$x_j = r - a_j\left(x_j^0 + \frac{a_i - r}{a_j} x_i^0\right) + \frac{a_i - r}{a_j} - \sum_{\substack{s=1\\ s\notin\{i,j\}}}^{n} a_s x_s^0 + \frac{r - a_i}{a_j} b = x_i^0$$

and $x_s = k_s = x_s^0,\ s \in \{1,...,n\} \setminus \{i,j\}$.

**Lemma 4.** Let $|r| \neq$ and $(i,j)$ be the pair of indices for which this minimum can be obtained. Again, let's consider the system of linear equations:

$$(3)\qquad \begin{cases} a_j t_h + r x_i + \sum_{\substack{s=1\\ s\notin\{i,j\}}}^{n} a_s x_s = b \\ t_h = x_j + \frac{a_i - r}{a_j} x_i \end{cases}$$

Then $x_e = x_e^0,\ e = \overline{1,n}$ is a particular integer solution for (2) if and only if $x_e = x_e^0$, $e \in \{1,...,n\} \setminus \{j\}$ and $t_h = t_h^0 = x_j^0 + \frac{a_i - r}{a_j} x_i$ is the particular integer solution of (3).

*Proof:*

$x_e = x_e^0,\ e = \overline{1,n}$ is a particular solution for (2) if and only if

$$\sum_{e=1}^{n} a_e x_e^0 = b \Leftrightarrow \sum_{\substack{s=1\\ s\notin\{i,j\}}}^{n} a_s x_s^0 + a_j\left(x_j^0 + \frac{a_i - r}{a_j} x_i^0\right) + r x_i^0 = b \Leftrightarrow$$

$$\Leftrightarrow a_j t_h^0 + r x_i^0 + \sum_{\substack{s=1\\ s\notin\{i,j\}}}^{n} a_s x_s^0 = b \qquad \text{and} \qquad t_h^0 = x_j^0 + \frac{a_i - r}{a_j} x_i^0 \in \mathbb{Z} \qquad \Leftrightarrow x_e = x_e^0,$$

$e \in \{1,...,n\} \setminus \{j\}$ and $t_h = t_h^0$ is a particular integer solution for (3).

**Lemma 5.** The previous algorithm is finite.

*Proof:*

When $|r| = 1$ the algorithm stops at step 3. We will discuss the case when $|r| \neq 1$. According to the definition of $r$, $|r| \in \mathbb{N}^*$. We will show that the row of $r-s$ successively obtained by following the algorithm several times is decreasing with cycle, and each cycle is not equal to the previous, by 1. Let $r_1$ be the first obtained by following the algorithm one time. $|r_1| \neq 1$ then go to step 4, and then step 5. According to lemma 2, $|r_1| < |a_i|$, $\forall i = \overline{1,n}$.

Now we shall follow the algorithm a second time, but this time for an equation in which $r_1$ (according to step 5) is substituted by $a_i$. Again, according to lemma 2, the new $|r|$ written $|r_2|$ will have the propriety: $|r_2| < |r_1|$. We will get to $|r| = 1$ because $|r| \geq 1$ and $|r| < \infty$, and if $|r| \neq 1$, following the algorithm once again we get $|r| < |r_1|$ and so on. Hence, the algorithm has a finite number of repetitions.

**Theorem of Correctness.** The previous algorithm calculates the general solution of the linear equation correctly (2).

*Proof:*

According to lemma 5 the algorithm is finite. From lemma 1 it follows that the set $M$ has a minimum, hence step 2 of the algorithm has meaning. When $|r| = 1$ it was shown in lemma 3 that step 3 of the algorithm calculates the general integer solution of the respective equation correctly the equation that appears at step 3). In lemma 4 it is shown that if $|r| \neq 1$ the substitutions steps 4 and 5 introduced in the initial equation, the general integer solution remains unchanged. That is, we pass from the initial equation to a linear system having the same general solution as the initial equation. The variable $h$ is a counter of the newly introduced variables, which are used to successively decompose the system in systems of two linear equations. The variable $p$ is a counter of the substitutions of variables (the relations, at a given moment between certain variables).

When the initial equation was decomposed to $|r| = 1$, we had to proceed in the reverse way, i.e. to compose its general integer solution. This reverse way is directed by the sub-steps 3(A), 3(B) and 3(C). The sub-step 3(D) has only an aesthetic role, i.e., to have the general solution under the form: $x_i = f_i(k_1, ..., k_{n-1})$, $i = \overline{1,n}$, $f_i$ being linear functions with integer number of coefficients. This "if possible" shows that substitutions are not always possible. But when they are we must make all possible substitutions.

**Note 1.** The previous algorithm can be easily introduced into a computer program.

**Note 2.** The previous algorithm is more "rapid" than that of W. Sierpinski's [1], i.e., the general integer solution is reached after a smaller number of iterations (or, at least, the same) for any linear equation (2).

In the first place, both methods aim at obtaining the coefficient $\pm 1$ for at least one unknown variable. While Sierpinski started only by chance, decomposing the greatest coefficient in the module (writing it under the form of a sum between a multiple of the following smaller coefficient (in the module) and the rest), in our algorithm this decomposition is not accidental but always seeks the smallest $|r|$ and also choose the coefficients $a_i$ and $a_j$ for which this minimum is achieved. That is, we test from the beginning the shortest way to the general integer solution. Sierpinski does not attempt to find the shortest way; he knows that his method will take him to the general integer solution of the equation and is not interested in how long it will take. However, when an algorithm is introduced into a computer program it is preferable that the process time should be as short as possible.

**Example 1.**

Let us solve in $\mathbb{Z}^3$ the equation $17x - 7y + 10z = -12$.

We apply the former algorithm.

1. $h = 1, p = 1$
2. $r = 3,\ i = 3,\ j = 2$
3. $|3| \neq 1$ go on to step 4.
4. (1) $y = t_1 - \dfrac{10-3}{-7} \cdot z = t_1 + z$
5. Assign

$y := t_1 \qquad h := 2$

$a_3 := 3 \qquad p := 2$

with the other coefficients and variables remaining unchanged, go back to step 2.

2. $r = -1,\ i = 1,\ j = 3$
3. $|-1| = 1$

$x = -1(-3t_2 - (-7t_1) - 12) = 3t_2 - 7t_1 - 12$

$$z = -1\left(17t_2 + (-7t_1) \cdot \frac{17-(-1)}{3} + \frac{-1-17}{3}(-12)\right) = 17t_2 + 42t_1 - 72$$

(A) We substitute the values of $x$ and $z$ thus determined into the only statement $(p)$ we have:

(1) $\quad y = t_1 + z = --17t_2 + 43t_1 - 72$

(B) The substitution is not possible.

(C) The substitution is not possible.

(D) The general integer solution of the equation is:

$$\begin{cases} x = 3k_1 - 7k_2 + 12 \\ y = -17k_1 + 43k_2 - 72 \\ z = -17k_1 + 42k_2 - 72; \qquad k_1, k_2 \in \mathbb{Z} \end{cases}$$

**REFERENCES:**

[1] Sierpinski, W, - Ce ştim şi ce nu ştim despre numerele prime? - Editura Stiinţifică, Bucharest, 1966.

[2] Creangă, I., Cazacu, C., Mihuţ, P., Opaiţ, Gh., Corina Reisher – Introducere în teoria numerelor, Ed. Did. şi Ped., Bucharest, 1965.

[3] Popovici, C. P. – Aritmetica şi teoria numerelor, Ed. Did. şi Ped., Bucharest, 1963.

[4] Smarandache, Florentin – Un algoritm de rezolvare în numere întregi a ecuaţiilor liniare.

# INTEGER NUMBER SOLUTIONS OF LINEAR SYSTEMS

## Definitions and Properties of the Integer Solution of a Linear System

Let's consider

$$(1) \qquad \sum_{j=1}^{n} a_{ij} x_j = b_i, \quad i = \overline{1,m}$$

a linear system with all coefficients being integer numbers (the case with rational coefficients is reduced to the same).

**Definition 1.** $x_j = x_j^0,\ j = \overline{1,n}$, is a particular integer solution of (1) if $x_j^0 \in \mathbb{Z},\ j = \overline{1,n}$ and $\sum_{j=1}^{n} a_{ij} x_j^0 = b_i,\ i = \overline{1,m}$.

Let's consider the functions $f_j : \mathbb{Z}^h \to \mathbb{Z},\ j = \overline{1,n}$, where $h \in \mathbb{N}^*$.

**Definition 2.** $x_j = f_j(k_1, ..., k_h),\ j = \overline{1,n}$, is the general integer solution for (1) if:

(a) $\sum_{j=1}^{n} a_{ij} f_j(k_1, ..., k_h) = b_i,\ i = \overline{1,m}$, irrespective of $(k_1, ..., k_h) \in \mathbb{Z}$;

(b) Irrespective of $x_j = x_j^0,\ j = \overline{1,n}$ a particular integer solution of (1) there is $(k_1^0, ..., k_h^0) \in \mathbb{Z}$ such that $f_j(k_1^0, ..., k_h^0) = x_j,\ j = \overline{1,n}$. (In other words the general solution that comprises all the other solutions.)

**Property 1.**

A general solution of a linear system of $m$ equations with $n$ unknowns, $r(A) = m < n$, is undetermined $(n-m)$-times.

*Proof:*

We assume by reduction ad absurdum that it is of order $r$, $1 \le r \le n-m$ (the case $r = 0$, i.e., when the solution is particular, is trivial). It follows that the general solution is of the form:

$$(S_1) \qquad \begin{cases} x_1 = u_{11} p_1 + ... + u_{1r} p_r + v_1 \\ \vdots \\ x_n = u_{n1} p_1 + ... + u_{nr} p_r + v_n, \quad u_{ih}, \forall i \in \mathbb{Z} \\ p_h = \text{parameters} \in \mathbb{Z} \end{cases}$$

We prove that the solution is undetermined $(n-m)$-times.

The homogeneous linear system (1), resolved in $r$ has the solution:

$$\begin{cases} x_1 = \dfrac{D_{m+1}^1}{D} x_{m+1} + \ldots + \dfrac{D_n^1}{D} x_n \\ \vdots \\ x_m = \dfrac{D_{m+1}^m}{D} x_{m+1} + \ldots + \dfrac{D_n^m}{D} x_n \end{cases}$$

Let $x_i = x_i^0,\ i = \overline{1,n}$, be a particular solution of the linear system (1).

Considering

$$\begin{cases} x_{m+1} = D \cdot k_{m+1} \\ \vdots \\ x_n = D \cdot k_n \end{cases}$$

we obtain the solution

$$\begin{cases} x_1 = D_{m+1}^1 \cdot k_{m+1} + \ldots + D_n^1 \cdot k_n + x_1^0 \\ \vdots \\ x_m = D_{m+1}^m \cdot k_{m+1} + \ldots + D_n^m \cdot k_n + x_m^0 \\ x_{m+1} = D \cdot k_{m+1} + x_{m+1}^0 \\ \vdots \\ x_n = D \cdot k_n + x_n^0, \qquad\qquad k_j = \text{parameters} \in \mathbb{Z} \end{cases}$$

which depends on the $n-m$ independent parameters, for the system (1). Let the solution be undetermined $(n-m)$-times:

(S2)
$$\begin{cases} x_1 = c_{11}k_1 + \ldots + c_{1n-m}k_{n-m} + d_1 \\ \vdots \\ x_n = c_{n1}k_1 + \ldots + u_{nn-m}k_{n-m} + d_n \\ c_{ij},\ d_i \in \mathbb{Z},\ k_j = \text{parameters} \in \mathbb{Z} \end{cases}$$

(There are such solutions, we have proved it before.) Let the system be:

$$\begin{cases} a_{11}x_1 + \ldots + a_{1n}x_n = b_1 \\ \vdots \\ a_{m1}x_1 + \ldots + a_{mn}x_n = b_m \end{cases}$$

$x_i$ = unknowns $\in \mathbb{Z}$, $a_{ij},\ b_i \in \mathbb{Z}$.

I. The case $b_i = 0,\ i = \overline{1,m}$ results in a homogenous linear system:

$a_{i1}x_i + \ldots + a_{in}x_n = 0;\ i = \overline{1,m}$.

(S2) $\Rightarrow a_{i1}\left(c_{i1}k_1 + \ldots + c_{1n-m}k_{n-m} + d_1\right) + \ldots + a_{in}\left(c_{n1}k_1 + \ldots + c_{nn-m}k_{n-m} + d_n\right) = 0$

$0 = \left(a_{i1}c_{11} + \ldots + a_{in}c_{n1}\right)k_1 + \ldots + \left(a_{i1}c_{1n-m} + \ldots + a_{in}c_{nn-m}\right)k_{n-m} + \left(a_{i1}d_1 + \ldots + a_{in}d_n\right)$

$\forall k_j \in \mathbb{Z}$

For $k_1 = \ldots = k_{n-m} = 0 \Rightarrow a_{i1}d_1 + \ldots + a_{in}d_n = 0$.

For $k_1 = ... = k_{h-1} = k_{h+1} = ... = k_{n-m} = 0$ and $k_h = 1 \Rightarrow$

$\Rightarrow \left(a_{i1}c_{ih} + ... + a_{in}c_{nh}\right) + \left(a_{i1}d_1 + ... + a_{in}d_d^{(n)}\right) = 0 \Rightarrow$

$a_{i1}c_{ih} + ... + a_{in}c_{nh} = 0,\ \forall\ i = \overline{1,m},\ \forall\ h = \overline{1,n-m}$.

The vectors

$$V_h = \begin{pmatrix} c_{1h} \\ : \\ : \\ c_{nh} \end{pmatrix},\ h = \overline{1,n-m}$$

are the particular solutions of the system.

$V_h,\ h = \overline{1,n-m}$ also linearly independent because the solution is undetermined $(n-m)$-times $\{V_1,...,V_{n-m}\} + d$ is a linear variety that includes the solutions of the system obtained from ($S_2$).

Similarly for ($S_1$) we deduce that

$$U_s = \begin{pmatrix} U_{1s} \\ : \\ : \\ U_{ns} \end{pmatrix},\ s = \overline{1,r}$$

are particular solutions of the given system and are linearly independent, because ($S_1$) is undetermined $(n-m)$-times, and $V = \begin{pmatrix} V_1 \\ : \\ : \\ V_n \end{pmatrix}$ is a solution of the given system.

**Case (a)** $U_1,...,U_r,\ V =$ linearly dependent, it follows that $\{U_1,...,U_r\}$ is a free sub-module of order $r < n - m$ of solutions of the given system, then, it follows that there are solutions that belong to $\{V_1,...,V_{n-m}\} + d$ and which do not belong to $\{U_1,...,U_r\}$, a fact which contradicts the assumption that ($S_1$) is the general solution.

**Case (b)** $U_1,...,U_r,\ V =$ linearly independent.

$\{U_1,...,U_r\}$+V is a linear variety that comprises the solutions of the given system, which were obtained from ($S_1$). It follows that the solution belongs to $\{V_1,...,V_{n-m}\} + d$ and does not belong to $\{U_1,...,U_r\}$+V, a fact which is a contradiction to the assumption that ($S_1$) is the general solution.

II. When there is an $i \in \overline{1,m}$ with $b_i \neq 0$ then non-homogeneous linear system

$$a_{i1}x_i + ... + a_{in}x_n = b_1,\ i = \overline{1,m}$$

($S_2$) $\Rightarrow a_{i1}\left(c_{11}k_1 + ... + c_{1n-m}k_{n-m} + d_1\right) + ... + a_{in}\left(c_{n1}k_1 + ... + c_{nn-m}k_{n-m} + d_n\right) = b_i$

it follows that

$$\Rightarrow \left(a_{i1}c_{11} + ... + a_{in}c_{n1}\right)k_1 + ... + \left(a_{i1}c_{1n-m} + ... + a_{in}c_{nn-m}\right)k_{n-m} + \left(a_{i1}d_1 + ... + a_{in}d_n\right) = b_i$$

For $k_1 = ... = k_{n-m} = 0 \Rightarrow a_{i1}d_1 + ... + a_{in}d_n = b_1$;
For $k_1 = ... = k_{j-1} = k_{j+1} = ... = k_{n-m} = 0$ and $k_j = 1 \Rightarrow$
$\Rightarrow \left(a_{i1}c_{1j} + ... + a_{in}c_{nj}\right) + \left(a_{in}d_1 + ... + a_{in}d_n\right) = b_i$ it follows that

$$\begin{cases} a_{i1}c_{1j} + ... + a_{in}c_{nj} = 0 \\ a_{i1}d_1 + ... + a_{in}d_n = b_i \end{cases}; \quad \forall i = \overline{1,m}, \ \forall j = \overline{1,n-m}.$$

$V_j = \begin{pmatrix} c_{1j} \\ : \\ c_{nj} \end{pmatrix}$, $j = \overline{1,n-m}$, are linearly independent because the solution ($S_2$) is undetermined $(n-m)$-times.

(?!) $V_j$, $j = \overline{1,n-m}$, and $d = \begin{pmatrix} d_1 \\ : \\ d_n \end{pmatrix}$

are linearly independent.

We assume that they are not linearly independent. It follows that

$$d = s_1V_1 + ... + s_{n-m}V_{n-m} = \begin{pmatrix} s_1c_{11} + ... + s_{n-m}c_{1n-m} \\ : \\ s_1c_{n1} + ... + s_{n-m}c_{nn-m} \end{pmatrix}.$$

Irrespective of $i = \overline{1,m}$:

$$b_1 = a_{i1}d_1 + ... + a_{in}d_n = a_{i1}\left(s_1c_{11} + ... + s_{n-m}c_{1n-m}\right) + ... + a_{in}\left(s_1c_{n1} + ... + s_{n-m}c_{nn-m}\right) =$$
$$= \left(a_{i1}c_{11} + ... + a_{in}c_{n1}\right)s_1 + .... + \left(a_{i1}c_{1n-m} + ... + a_{in}c_{nn-m}\right)s_{n-m} = 0.$$

Then, $b_i = 0$, irrespective of $i = \overline{1,m}$, contradicts the hypothesis (that there is an $i \in \overline{1,m}$, $b_i \neq 0$). It follows that $V_1,...,V_{n-m}, d$ are linearly independent.

$\{V_1,...,V_{n-m}\} + d$ is a linear variety that contains the solutions of the non-homogeneous system, solutions obtained from ($S_2$). Similarly it follows that $\{G_1,...,G_r\} + V$ is a linear variety containing the solutions of the non-homogeneous system, obtained from ($S_1$).

$n - m > r$ it follows that there are solutions of the system that belong to

---

"?!" means "to prove that"

$\{V_1,...,V_{n-m}\} + d$ and which do not belong to $\{G_1,...,G_r\} + V$, this contradicts the fact that ($S_1$) is the general solution. Then, it shows that the general solution depends on the $n - m$ independent parameters.

**Theorem 1.** The general solution of a non-homogeneous linear system is equal to the general solution of an associated linear system plus a particular solution of the non-homogeneous system.

*Proof:*

Let's consider the homogeneous linear solution:

$$\begin{cases} a_{11}x_1 + ... + a_{1n}x_n = 0 \\ : \\ a_{m1}x_1 + ... + a_{mn}x_n = 0 \end{cases}, \quad (AX = 0)$$

with the general solution:

$$\begin{cases} x_1 = c_{11}k_1 + ... + c_{1n-m}k_{n-m} + d_1 \\ : \\ x_n = c_{n1}k_1 + ... + c_{nn-m}k_{n-m} + d_n \end{cases}$$

and

$$\begin{cases} x_1 = x_1^0 \\ : \\ x_n = x_n^0 \end{cases}$$

with the general solution a particular solution of the non-homogeneous linear system $AX = b$;

$$(?!) \qquad \begin{cases} x_1 = c_{11}k_1 + ... + c_{1n-m}k_{n-m} + d + x_1^0 \\ : \\ x_n = c_{n1}k_1 + ... + c_{nn-m}k_{n-m} + d_n + x_n^0 \end{cases}$$

is a solution of the non-homogeneous linear system.

We note:

$$A = \begin{pmatrix} a_{11} ... a_{1n} \\ : \\ a_{m1} ... a_{mn} \end{pmatrix}, \; X = \begin{pmatrix} x_1 \\ : \\ x_n \end{pmatrix}, \; b = \begin{pmatrix} b_1 \\ : \\ b_m \end{pmatrix}, \; 0 = \begin{pmatrix} 0 \\ : \\ 0 \end{pmatrix}$$

(vector of dimension $m$),

$$K = \begin{pmatrix} k_1 \\ : \\ k_{n-m} \end{pmatrix}, \; C = \begin{pmatrix} c_{11} ... c_{1n-m} \\ : \\ c_{n1} ... c_{nn-m} \end{pmatrix}, \; d = \begin{pmatrix} d_1 \\ : \\ d_n \end{pmatrix}, \; x^0 = \begin{pmatrix} x_1^0 \\ : \\ x_n^0 \end{pmatrix};$$

$$AX = A\left(Ck + d + x^0\right) = A\left(Ck + d\right) + AX^0 = b + 0 = b$$

We will prove that irrespective of

$$x_1 = y_1^0$$
$$:$$
$$x_n = y_n^0$$

there is a particular solution of the non-homogeneous system

$$\begin{cases} k_1 = k_1^0 \in \mathbb{Z} \\ \vdots \\ k_{n-m} = k_{n-m}^0 \in \mathbb{Z} \end{cases},$$

with the property:

$$\begin{cases} x_1 = c_{11}k_1^0 + \ldots + c_{1n}k_{n-m}^0 + d_1 + x_1^0 = y_1^0 \\ \vdots \\ x_n = c_{n1}k_1^0 + \ldots + c_{nn-m}k_{n-m}^0 + d_1 + x_n^0 = y_n^0 \end{cases}$$

We note $Y^0 = \begin{pmatrix} y_1^0 \\ \vdots \\ y_n^0 \end{pmatrix}$.

We'll prove that those $k_j^0 \in \mathbb{Z}, j = \overline{1, n-m}$ are those for which $A\left(CX^0 + d\right) = 0$ (there are such $X_j^0 \in \mathbb{Z}$ because

$$\begin{cases} x_1 = 0 \\ \vdots \\ x_n = 0 \end{cases}$$

is a particular solution of the homogeneous linear system and $X = CK + d$ is a general solution of the non-homogeneous linear system)

$$A\left(CK^0 + d + X^0 - Y^0\right) = A\left(CK^0 + d\right) + AX^0 - AY^0 = 0 + b - b = 0 \quad .$$

**Property 2** The general solution of the homogeneous linear system can be written under the form:

(SG)

(2) $$\begin{cases} x_1 = c_{11}k_1 + \ldots + c_{1n-m}k_{n-m} \\ \vdots \\ x_n = c_{n1}k_1 + \ldots + c_{nn-m}k_{n-m} \end{cases}$$

$k_j$ is a parameter that belongs to $\mathbb{Z}$ (with $d_1 = d_2 = \ldots = d_n = 0$).

*Poof:*

(SG) = general solution. It results that (SG) is undetermined $(n-m)$-times.

Let's consider that (SG) is of the form

(3) $$\begin{cases} x_1 = c_{11}p_1 + \ldots + c_{1n-m}p_{n-m} + d_1 \\ \vdots \\ x_n = c_{n1}p_1 + \ldots + c_{nn-m}p_{n-m} + d_n \end{cases}$$

with not all $d_i = 0$; we'll prove that it can be written under the form (2); the system has the trivial solution

$$\begin{cases} x_1 = 0 \in \mathbb{Z} \\ \vdots \\ x_n = 0 \in \mathbb{Z} \end{cases};$$

it results that there are $p_j \in \mathbb{Z},\ j = \overline{1, n-m}$,

$$(4) \qquad \begin{cases} x_1 = c_{11} p_1^0 + \ldots + c_{1n-m} p_{n-m}^0 + d_1 = 0 \\ \vdots \\ x_n = c_{n1} p_1^0 + \ldots + c_{nn-m} p_{n-m}^0 + d_n = 0 \end{cases}$$

Substituting $p_j = k_j + p_j^0,\ j = \overline{1, n-m}$ in (3)

$$\left.\begin{matrix} k_j \in \mathbb{Z} \\ p_j^0 \in \mathbb{Z} \end{matrix}\right\} \Rightarrow p_j \in \mathbb{Z},$$

$$\left.\begin{matrix} p_j \in \mathbb{Z} \\ p_j^0 \in \mathbb{Z} \end{matrix}\right\} \Rightarrow k_j = p_j - p_j^0 \in \mathbb{Z}$$

which means that that they do not make any restrictions.

It results that

$$\begin{cases} x_1 = c_{11} k_1 + \ldots + c_{1n-m} k_{n-m} + \left(c_{11} p_1^0 + \ldots + c_{1n-m} p_{n-m}^0 + d_1\right) \\ \vdots \\ x_n = c_{n1} k_1 + \ldots + c_{nn-m} k_{n-m} + \left(c_{n1} p_1^0 + \ldots + c_{nn-m} p_{n-m}^0 + d_n\right) \end{cases}$$

But

$$c_{h1} p_1^0 + \ldots + c_{hn-m} p_{n-m}^0 + d_h = 0,\ h = \overline{1, n} \text{ (from (4)).}$$

Then the general solution is of the form:

$$\begin{cases} x_1 = c_{11} k_1 + \ldots + c_{1n-m} k_{n-m} \\ \vdots \\ x_n = c_{n1} k_1 + \ldots + c_{nn-m} k_{n-m} \end{cases}$$

$k_j$ = parameters $\in \mathbb{Z}$, $j = \overline{1, n-m}$; it results that $d_1 = d_2 = \ldots = d_n = 0$.

**Theorem 2.** Let's consider the homogeneous linear system:

$$\begin{cases} a_{11} x_1 + \ldots + a_{1n} x_n = 0 \\ \vdots \\ a_{m1} x_1 + \ldots + a_{mn} x_n = 0 \end{cases},$$

$r(A) = m$, $(a_{h1}, \ldots, a_{hn}) = 1$, $h = \overline{1, m}$ and the general solution

$$\begin{cases} x_1 = c_{11} k_1 + \ldots + c_{1n-m} k_{n-m} \\ \vdots \\ x_n = c_{n1} k_1 + \ldots + c_{nn-m} k_{n-m} \end{cases}$$

then

$$\left(a_{h1},...,a_{hi-1},a_{hi+1},...,a_{hn}\right)\Big|\left(c_{i1},...,c_{in-m}\right)$$

irrespective of $h=\overline{1,m}$ and $i=\overline{1,n}$.

*Proof:*

Let's consider some arbitrary $h\in\overline{1,m}$ and some arbitrary $i\in\overline{1,n}$;

$$a_{h1}x_1+...+a_{hi-1}x_{i-1}+a_{hi+1}x_{i+1}+...+a_{hn}x_n=a_{hi}x_i.$$

Because

$$\left(a_{h1},...,a_{hi-1},a_{hi+1},...,a_{hn}\right)\Big|a_{hi}$$

it results that

$$d=\left(a_{h1},...,a_{hi-1},a_{hi+1},...,a_{hn}\right)\Big|x_i$$

irrespective of the value of $x_i$ in the vector of particular solutions.

For $k_2=k_3=...=k_{n-m}=0$ and $k_1=1$ we obtain the particular solution:

$$\begin{cases}x_1=c_{11}\\ \vdots\\ x_i=c_{i1}\quad\Rightarrow d\mid c_{i1}\\ \vdots\\ x_n=c_{n1}\end{cases}$$

For $k_1=k_2=...=k_{n-m-1}=0$ and $k_{n-m}=1$ it results the following particular solution:

$$\begin{cases}x_1=c_{1n-m}\\ \vdots\\ x_i=c_{in-m}\quad\Rightarrow d\mid c_{in-m};\\ \vdots\\ x_n=c_{nn-m}\end{cases}$$

hence

$$d\mid c_{ij},\ j=\overline{1,n-m}\Rightarrow d\big|\left(c_{i1},...,c_{in-m}\right).$$

**Theorem 3.**

If

$$\begin{cases}x_1=c_{11}k_1+...+c_{1n-m}k_{n-m}\\ \vdots\\ x_n=c_{n1}k_1+...+c_{nn-m}k_{n-m}\end{cases}$$

$k_j$ = parameters $\in\mathbb{Z}$, $c_{ij}\in\mathbb{Z}$ being given, is the general solution of the homogeneous linear system

$$\begin{cases}a_{11}x_1+...+a_{1n}x_n=0\\ \vdots\\ a_{m1}x_1+...+a_{mn}x_n=0\end{cases},\quad r(A)=m<n$$

then $\left(c_{1j},...,c_{nj}\right)=1$, $\forall j=\overline{1,n-m}$.

*Proof:*

We assume, by reduction ad absurdum, that there is $j_0 \in \overline{1,n-m}:\left(c_{1j_0},...,c_{nj_0}\right)=d$ we consider the maximal co-divisor $>0$; we reduce to the case when the maximal co-divisor is $-d$ to the case when it is equal to $d$ (non restrictive hypothesis); then the general solution can be written under the form:

$$(5)\qquad \begin{cases} x_1 = c_{11}k_1 + ... + c'_{1j_0}dk_{j_0} + ... + c_{1n-m}k_{n-m} \\ : \\ x_n = c_{n1}k_1 + ... + c'_{nj_0}dk_{j_0} + ... + c_{nn-m}k_{n-m} \end{cases}$$

where $d=\left(c_{ij_0},..,c_{nj_0}\right)$, $c_{ij_0}=d\cdot c'_{ij_0}$ and $\left(c'_{ij_0},...,c'_{nj_0}\right)=1$.

We prove that

$$\begin{cases} x_1 = c'_{1j_0} \\ : \\ x_n = c'_{nj_0} \end{cases}$$

is a particular solution of the homogeneous linear system.

We'll note:

$$C=\begin{pmatrix} c_{11} \;...\; c'_{ij_0} \;\; d \;...\; c_{1n-m} \\ : \qquad\quad : \qquad\qquad : \\ c_{n1} \;...\; c'_{nj_0} \;\; d \;...\; c_{nn-m} \end{pmatrix},\; k=\begin{pmatrix} k_1 \\ : \\ k_{j_0} \\ : \\ k_{n-m} \end{pmatrix}$$

$x=C\cdot k$ the general solution.

We know that $AX=0 \Rightarrow A(CK)=0$, $A=\begin{pmatrix} a_{11} ... \; a_{1n} \\ : \\ a_{n1} ... \; a_{mn} \end{pmatrix}$.

We assume that the principal variables are $x_1,...,x_m$ (if not, we have to renumber). It follows that $x_{m+1},...,x_n$ are the secondary variables.

For $k_1=...=k_{j_0-1}=k_{j_0+1}=...=k_{n-m}=0$ and $k_{j_0}=1$ we obtain a particular solution of the system

$$\begin{cases} x_1 = c'_{1j_0}d \\ : \\ x_n = c'_{nj_0}d \end{cases} \Rightarrow 0 = A\begin{pmatrix} c'_{1j_0}d \\ : \\ c'_{nj_0}d \end{pmatrix} = d\cdot A\begin{pmatrix} c'_{1j_0} \\ : \\ c'_{nj_0} \end{pmatrix} \Rightarrow A\begin{pmatrix} c'_{1j_0} \\ : \\ c'_{nj_0} \end{pmatrix} = 0 \Rightarrow \begin{cases} x_1 = c'_{1j_0} \\ : \\ x_n = c'_{nj_0} \end{cases}$$

is the particular solution of the system.

We'll prove that this particular solution cannot be obtained by

$$(6)\qquad \begin{cases} x_1 = c_{11}k_1 + ... + c'_{1j_0}dk_{j_0} + ... + c_{1n-m}k_{n-m} = c'_{1j_0} \\ : \\ x_n = c_{n1}k_1 + ... + c'_{nj_0}dk_{j_0} + ... + c_{nn-m}k_{n-m} = c'_{nj_0} \end{cases}$$

$$(7)\qquad \begin{cases} x_{m+1} = c_{m+1}k_1 + ... + c'_{m+1}dk_{j_0} + ... + c_{m+1,n-m}k_{n-m} = c'_{m+1\,j_0} \\ : \\ x_n = c_{n1}k_1 + ... + c'_{nj_0}dk_{j_0} + ... + c_{nn-m}k_{n-m} = c'_{nj_0} \end{cases}$$

$$\Rightarrow\ k_{j_0} = \frac{\begin{vmatrix} c_{m+1,1} & ... & c_{m+1\,j} & ... & c_{m+1,n-m} \\ : & & : & 0. & : \\ c_{h,1} & ... & c_{nj} & ... & c_{n,n-m} \end{vmatrix}}{\begin{vmatrix} c_{m+1,1} & ... & c'_{m+1\,j_0}d & ... & c_{m+1,n-m} \\ : & & : & 0. & : \\ c_{h,1} & ... & c'_{nj}d & ... & c_{n,n-m} \end{vmatrix}} = \frac{1}{d} \notin \mathbb{Z}$$

(because $d \neq 1$).

It is important to point out the fact that those $k_j = k_j^0$, $j = \overline{1, n-m}$, that satisfy the system (7) also satisfy the system (6), because, otherwise (6) would not satisfy the definition of the solution of a linear system of equations (i.e., considering the system (7) the hypothesis was not restrictive). From $X_{j_0} \in \mathbb{Z}$ follows that (6) is not the general solution of the homogeneous linear system contrary to the hypothesis); then $(c_{1j}, ..., c_{nj}) = 1$, irrespective of $j = \overline{1, n-m}$.

**Property 3.** Let's consider the linear system

$$\begin{cases} a_{11}x_1 + ... + a_{1n}x_n = b_1 \\ : \\ a_{m1}x_1 + ... + a_{mn}x_n = b_m \end{cases}$$

$a_{ij}, b_i \in \mathbb{Z}$, $r(A) = m < n$, $x_j =$ unknowns $\in \mathbb{Z}$

Resolved in $\mathbb{R}$, we obtain

$$\begin{cases} x_1 = f_1(x_{m+1}, ..., x_n) \\ : \\ x_m = f_m(x_{m+1}, ..., x_n) \end{cases}, \quad x_1, ..., x_m \text{ are the main variables,}$$

where $f_i$ are linear functions of the form:

$$f_i = \frac{c^i_{m+1}x_{m+1} + ... + c^i_n x_n + e_i}{d_i},$$

where $c^i_{m+j}$, $d_i$, $e_i \in \mathbb{Z}$; $i = \overline{1, m}$, $j = \overline{1, n-m}$.

If $\frac{e_i}{d_i} \in \mathbb{Z}$ irrespective of $i = \overline{1,m}$ then the linear system has integer solution.

*Proof:*

For $1 \le i \le m,\ x_i \in \mathbb{Z}$, then $f_j \in \mathbb{Z}$. Let's consider

$$\begin{cases} x_{m+1} = u_{m+1}k_{m+1} \\ \vdots \\ x_n = u_n k_n \\ \vdots \\ x_1 = v^1_{m+1}k_{m+1} + ... + v^1_n k_n + \frac{e_1}{d_1} \\ \vdots \\ x_m = v^m_{m+1}k_{m+1} + ... + v^m_n k_n + \frac{e_m}{d_m} \end{cases}$$

a solution, where $u_{m+1}$ is the maximal co-divisor of the denominators of the fractions $\frac{c^i_{m+j}}{d_i},\ i = \overline{1,m},\ j = \overline{1,n-m}$ calculated after their complete simplification.

$v^i_{m+j} = \frac{c^i_{m+j}u_{m+j}}{d_i} \in \mathbb{Z}$ is a $(n-m)$-times undetermined solution which depends on $n-m$ independent parameters $\left(k_{m+1},...,k_n\right)$ but is not a general solution.

**Property 4.** Under the conditions of property 3, if there is an $i_0 \in \overline{1,m} : f_{i_0} = u^{i_0}_{m+1}x_{m+1} + ... + u^{i_0}_n x_n + \frac{e_{i_0}}{d_{i_0}}$ with $u^{i_0}_{m+j} \in \mathbb{Z},\ j = \overline{1,n-m}$, and $\frac{e_{i_0}}{d_{i_0}} \notin \mathbb{Z}$ then the system does not have integer solution.

*Proof:*

$\forall x_{m+1},...,x_n$ in $\mathbb{Z}$, it results that $x_{i_0} \notin \mathbb{Z}$.

**Theorem 4.** Let's consider the linear system

$$\begin{cases} a_{11}x_1 + ... + a_{1n}x_n = b_1 \\ \vdots \\ a_{m1}x_1 + ... + a_{mn}x_n = b_m \end{cases}$$

$a_{ij}, b_i \in \mathbb{Z},\ x_j =$ unknowns $\in \mathbb{Z},\ r(A) = m < n$. If there are indices $1 \le i_1 < ... < i_m \le n$, $i_h \in \{1,2,..,n\},\ h = \overline{1,m}$, with the property:

$$\Delta = \begin{vmatrix} a_{1i_1} & \dots & a_{1i_m} \\ \vdots & & \vdots \\ a_{mi_1} & \dots & a_{mi_m} \end{vmatrix} \neq 0 \text{ and}$$

$$\Delta_{x_{i_1}} = \begin{vmatrix} b_1 & a_{1i_2} & \dots & a_{1i_m} \\ \vdots & \vdots & & \vdots \\ b_m & a_{mi_2} & \dots & a_{mi_m} \end{vmatrix} \text{ is divided by } \Delta$$

.
.
.

$$\Delta_{x_{i_m}} = \begin{vmatrix} a_{1i_1} & \dots & a_{1i_{m-1}} & b_1 \\ \vdots & & \vdots & \vdots \\ a_{mi_1} & \dots & a_{mi_{m-1}} & b_m \end{vmatrix} \text{ is divided by } \Delta$$

then the system has integer number solutions.

*Proof:*

We use property 3

$$d_i = \Delta,\ i = \overline{1,m};\ e_{i_h} = \Delta_{x_{i_h}},\ h = \overline{1,m}$$

**Note 1.** It is not true in the reverse case.

**Consequence 1.** Any homogeneous linear system has integer number solutions (besides the trivial one); $r(A) = m < n$.

*Proof:*

$$\Delta_{x_{i_h}} = 0 : \Delta, \text{ irrespective of } h = \overline{1,m}.$$

**Consequence 2.** If $\Delta = \pm 1$, it follows that the linear system has integer number solutions.

*Proof:*

$\Delta_{x_{i_h}} : (\pm 1)$, irrespective of $h = \overline{1,m}$;

$\Delta_{x_{i_h}} \in \mathbb{Z}$.

# FIVE INTEGER NUMBER ALGORITHMS TO SOLVE LINEAR SYSTEMS

This section further extends the results obtained in chapters 4 and 5 (from linear equation to linear systems). Each algorithm is thoroughly proved and then an example is given.

Five integer number algorithms to solve linear systems are further given.

**Algorithm 1**. (Method of Substitution)

(Although simple, this algorithm requires complex computations but is, nevertheless, easy to implement into a computer program).

Some integer number equation are initially solved (which is usually simpler) by means of one of the algorithms 4 or 5. (If there is an equation of the system which does not have integer systems, then the integer system does not have integer systems, then Stop.) The general integer solution of the equation will depend on $n-1$ integer number parameters (see [5]):

$$(p_1) \qquad x_{i_1} = f_{i_1}^{(1)}\left(k_1^{(1)},...,k_{n-1}^{(1)}\right),\ i=\overline{1,n},$$

where all functions $f_{i_1}^{(1)}$ are linear and have integer number coefficients.

This general integer number system $(p_1)$ is introduced into the other $m-1$ equations of the system. We obtain a new system of $m-1$ equations with $n-1$ unknown variables:

$$k_{i_1}^{(1)},\ i_1=\overline{1,n-1},$$

which is also to be solved with integer numbers. The procedure is similar. Solving a new equation, we obtain its general integer solution:

$$(p_2) \qquad k_{i_2}^{(1)} = f_{i_2}^{(2)}\left(k_1^{(2)},...,k_{n-2}^{(2)}\right),\ i_2=\overline{1,n-1},$$

where all functions $f_{i_2}^{(2)}$ are linear, with integer number coefficients. (If, along this algorithm we come across an equation, which does not have integer solutions, then the initial system does not have integer solution. Stop.)

In the case that all solved equations had integer system at step $(j)$, $1 \le j \le r$ ($r$ being of the same rank as the matrix associated to the system) then:

$$(p_j) \qquad k_{i_j}^{(j-1)} = f_{i_j}^{(j)}\left(k_1^{(j)},...,k_{n-j}^{(j)}\right),\ i_j=\overline{1,n-j+1},$$

$f_{i_j}^{(j)}$ are linear functions and have integer number coefficients.

Finally, after $r$ steps, and if all solved equations had integer solutions, we obtain the integer solution of one equation with $n-r+1$ unknown variables.

The system will have integer solutions if and only if in this last equation has integer solutions.

If it does, let its general integer solution be:

$$(p_r) \qquad k_{i_r}^{(r-1)} = f_{i_r}^{(r)}\left(k_1^{(r)},...,k_{n-1}^{(r)}\right),\ i_r=\overline{1,n-r+1},$$

where all $f_{i_r}^{(r)}$ are linear functions with integer number coefficients.

We'll present now the reverse procedure as follows.

We introduce the values of $k_{i_r}^{(r-1)}$, $i_r = \overline{1, n-r+1}$, at step $p_r$ in the values of

$$k_{i_{r-1}}^{(r-2)},\ i_{r-1} = \overline{1, n-r+2}$$

from step $(p_{r-1})$.

It follows:

$$k_{i_{r-1}}^{(r-2)} = f_{i_{r-1}}^{(r-1)}\left(f_1^{(r)}\left(k_1^{(r)},...,k_{n-r}^{(r)}\right),...,f_{n-r+1}^{(r)}\left(k_1^{(r)},...,k_{n-r}^{(r)}\right)\right) = g_{i_{r-1}}^{(r-1)}\left(k_1^{(r)},...,k_{n-r}^{(r)}\right),$$

$i_{r-1} = \overline{1, n-r-1}$, from which it follows that $g_{i_r}^{(r-1)}$ are linear functions with integer number coefficients.

Then follows those $(p_{r-2})$ in $(p_{r-e})$ and so on, until we introduce the values obtained at step $(p_2)$ in those from the step $(p_1)$.

It will follow:

$$x_{i_j} = g_i^1\left(k_1^{(r)},...,k_{n-r}^{(r)}\right)$$

notation $g_{i_1}\left(k_1,...,k_{n-r}\right)$, $i = \overline{1, n}$, with all $g_{i_1}$ most obviously, linear functions with integer number coefficients (the notation was made for simplicity and aesthetical aspects). This is, thus, the general integer solution, of the initial system.

**The correctness of Algorithm 1**.

The algorithm is finite because it has $r$ steps on the forward way and $r-1$ steps on the reverse, $(r < +\infty)$. Obviously, if one equation of one system does not have (integer number) solutions then the system does not have solutions either.

Writing $S$ for the initial system and $S_j$ the system resulted from step $(p_j)$, $1 \le j \le r-2$, it follows that passing from $(p_j)$ to $(p_{j+1})$ we pass from a system $S_j$ to a system $S_{j+1}$ equivalent from the point of view of the integer number solution, i.e.

$$k_{i_j}^{(j-1)} = t_{i_j}^0,\ \ i_j = \overline{1, n-j+1},$$

which is a particular integer solution of the system $S_j$ if and only if

$$k_{i_{j+1}}^{(j)} = h_{i_{j+1}}^0,\ i_{j+1} = \overline{1, n-j},$$

is a particular integer solution of the system $S_{j+1}$ where

$$k_{i_{j+1}}^0 = f_{i_{j+1}}^{(j+1)}\left(t_1^0,...,t_{n-j+1}^0\right),\ i_{j+1} = \overline{1, n-j}.$$

Hence, their general integer solutions are also equivalent (considering these substitutions). Such that, in the end, resolving the initial system $S$ is equivalent with solving the equation (of the system consisting of one equation) $S_{r-1}$ with integer number coefficients. It follows that the system $S$ has integer number solution if and only if the systems $S_j$ have integer number solution, $1 \le j \le r-1$.

**Example 1.** By means of algorithm 1, let us calculate the integer number solution of the following system:

$$(S) \qquad \begin{cases} 5x - 7y - 2z + 6w = 6 \\ -4x + 6y - 3z + 11w = 0 \end{cases}$$

Solution: We solve the first integer number equation. We obtain the general solution (see [4] or [5]):

$$(p_1) \qquad \begin{cases} x = t_1 + 2t_2 \\ y = t_1 \\ z = -t_1 + 5t_2 + 3t_3 - 3 \\ w = t_3 \end{cases}$$

where $t_1,\ t_2,\ t_3 \in \mathbb{Z}$.

Substituting in the second, we obtain the system:

$(S_1)$ $\qquad 5t_1 - 23t_2 + 2t_3 + 9 = 0$.

Solving this integer equation we obtain its general integer solution:

$$(p_2) \qquad \begin{cases} t_1 = k_1 \\ t_2 = k_1 + 2k_2 + 1 \\ t_3 = 9k_1 + 23k_2 + 7 \end{cases}$$

where $k_1,\ k_2 \in \mathbb{Z}$.

The reverse way. Substituting $(p_2)$ in $(p_1)$ we obtain:

$$\begin{cases} x = 3k_1 + 4k_2 + 2 \\ y = k_1 \\ z = 31k_1 + 79k_2 + 23 \\ w = 9k_1 + 23k_2 + 7 \end{cases}$$

where $k_1,\ k_2 \in \mathbb{Z}$, which is the general integer solution of the initial system $(S)$. Stop.

**Algorithm 2.**

**Input**

A linear system (1) without all $a_{ij} = 0$.

**Output**

We decide on the possibility of an integer solution of this system. If it is possible, we obtain its general integer solution.

**Method**

1. $t = 1,\ h = 1,\ p = 1$

2. (A) Divide each equation by the largest co-divisor of the coefficients of the unknown variables. If you do not obtain an integer quotient for at least one equation, then the system does not have integer solutions. Stop.

(B) If there is an inequality in the system, then the system does not have integer solutions. Stop.

(C) If repetition of more equations occurs, keep one and if an equation is an identity, remove it from the system.

3. If there is $(i_0, j_0)$ such that $\left|a_{i_0 j_0}\right| = 1$ then obtain the value of the variable $x_{j_0}$ from equation $i_0$; statement $(T_t)$.

Substitute this statement (where possible) in the other equations of the system and in the statement $(T_{t-1})$, $(H_h)$ and $(P_p)$ for all $i$, $h$, and $p$. Consider $t := t+1$, remove equation $(i_0)$ from the system. If there is no such a pair, go to step 5.

4. Does the system (left) have at least one unknown variable? If it does, consider the new data and go on to step 2. If it does not, write the general integer solution of the system substituting $k_1, k_2, \ldots$ for all variables from the right term of each expression which gives the value of the unknowns of the initial system. Stop.

5. Calculate

$$a = \min_{i, j_1, j_2} \left\{ |r| \, a_{ij_1} \equiv r \left( \operatorname{mod} a_{ij_2} \right), \; 0 < |r| < \left| a_{ij_2} \right| \right\}$$

and determine the indices $i, j_1, j_2$ as well as the $r$ for which this minimum can be calculated. (If there are more variables, choose one arbitrarily.)

6. Write: $x_{j_2} = t_h \dfrac{a_{ij_1} - r}{a_{ij_2}} x_{ij_2}$, statement $(H_h)$. Substitute this statement (where possible in all the equations of the system and in the statements $(T_t)$, $(H_h)$ and $(P_p)$ for all $t$, $h$, and $p$.

7. (A) If $a \neq 1$, consider $x_{j_2} := t_h$, $h := h+1$, and go on to step 2.

(B) If $a = 1$, then obtain the value of $x_{j_1}$ from the equation $(i)$; statement $(P_p)$. Substitute this statement (where possible in the other equations of the system and in the relations $(T_t)$, $(H_h)$ and $(P_{p-1})$ for all $t$, $h$, and $p$.

Remove the equation $(i)$ from the system.

Consider $h := h+1$, $p := p+1$, and go back to step 4.

**The correctness of algorithm 2.** Let consider system (1).

**Lemma 1.** We consider the algorithm at step 5. Also, let

$$M = \left\{ |r|, \; a_{ij_1} \equiv r \left( \operatorname{mod} a_{ij_2} \right), \; 0 < |r| < \left| a_{ij_2} \right|, \; i, j_1, j_2 = 1, 2, 3, \ldots \right\}.$$

Then $M \neq \emptyset$.

*Proof***:**

Obviously, $M$ is finite and $M \subset \mathbb{N}^*$. Then, $M$ has a minimum if and only if $M \neq \emptyset$. We suppose, conversely, that $M = \emptyset$. Then

$$a_{ij_2} \equiv 0 (\operatorname{mod} a_{ij_2}), \; \forall \; i, j_1, j_2.$$

It follows as well that

$$a_{ij_2} \equiv 0 (\operatorname{mod} a_{ij_1}), \; \forall \; i, j_1, j_2.$$

That is

$$\left| a_{ij_1} \right| = \left| a_{ij_2} \right|, \; \forall \; i, j_1, j_2.$$

We consider an $i_0$ arbitrary but fixed. It is clear that

$$(a_{i_0 1},...,a_{i_0 n}) : \ a_{i_0 j} \neq 0, \ \forall j$$

(because the algorithm has passed through the sub-steps 2(B) and 2(C). But, because it has also passed through step 3, it follows that

$$\left|a_{i_0 j}\right| \neq 1, \ \forall j ,$$

but as it previously passed through step 2(A), it would result that

$$\left|a_{i_0 j}\right| = 1, \ \forall j .$$

This contradiction shows that the assumption is false.

**Lemma 2.** Let's consider $a_{i_0 j_1} \equiv r(\mathrm{mod}\, a_{ij_2})$. Substitute

$$x_{j_2} = t_h - \frac{a_{i_0 j} - r}{a_{i_0 j_2}} x_{j_1}$$

in system (A) obtaining system (B). Then

$$x_j = x_j^0, \ j = \overline{1,n}$$

is the particular integer solution of (A) if and only if

$$x_j = x_j^0, \ j \neq j_2 \ \text{ and } \ t_h = x_{j_2}^0 - \frac{a_{i_0 j_1} - r}{a_{i_0 j_2}}$$

is the particular integer solution of (B).

**Lemma 3.** Let $a_1 \neq$ and $a_2$ obtained at step 5.

Then $0 < a_2 < a_1$

*Proof:*

It is sufficient to show that $a_1 < \left|a_{ij}\right|, \ \forall \ i, j$ because in order to get $a_2$, step 6 is obligatory, when the coefficients if the new system are calculated, $a_1$ being equal to a coefficient form the new system (equality of modules), the coefficient on $(i_0 \ j_1)$.

Let $a_{i_0 j_0}$ with the property $\left|a_{i_0 j_0}\right| \leq a_1$.

Hence, $a_1 \geq \left|a_{i_0 j}\right| = \min\left\{\left|a_{i_0 j}\right|\right\}$. Let $a_{i_0 j_s}$ with $\left|a_{i_0 j_s}\right| > \left|a_{ij_m}\right|$; there is such an element because $\left|a_{i_0 j_m}\right|$ is the minimum of the coefficients in the module and not all $\left|a_{i_0 j}\right|, \ j = \overline{1,n}$ are equal (conversely, it would result that $(a_{i_0 j},....,a_{i_0 n}) \sim a_{i_0 j}, \ \forall j \in \overline{1,r}$, the algorithm passed through sub-step 2(A) has simplified each equation by the maximal co-divisor of its coefficients; hence, it would follow that $\left|a_{i_0 j}\right| = 1, \ \forall j = \overline{1,n}$, which, again, cannot be real because the algorithm also passed through step 3). Out of the coefficients $a_{i_0 j_m}$ we choose one with the property $a_{i_0 j_{s_0}} \neq M a_{i_0 j_m}$ there is such an element (contrary, it would result $(a_{i_0 j},...,a_{i_0 n}) \sim \left|a_{i_0 j_m}\right|$ but the algorithm has also passed through step 2(A) and it would mean that $\left|a_{i_0 j_m}\right| = 1$ which contradicts step 3 through which the algorithm has also passed).

Considering $q_0 = \left[ a_{i_0 j_{s_0}} / a_{i_0 j_m} \right] \in \mathbb{Z}$ and $r = a_{i_0 j_{s_0}} - q_0 a_{i_0 j_m} \in \mathbb{Z}$, we have $a_{i_0 j_{s_0}} \equiv r_0 (\text{mod}\, a_{i_0 j_m})$ and $0 < |r_0| < \left|a_{i_0 j_m}\right| < \left|a_{i_0 j_0}\right| \le a_1$. We have, thus, obtained an $r_0$ with $|r_0| < a_1$, which is in contradiction with the very definition of $a_1$. Thus $a_1 < |a_{ij}|$, $\forall i, j$.

**Lemma 4.** Algorithm 2 is finite.

*Proof:*

The functioning of the algorithm is meant to transform a linear system of $m$ equations and $n$ unknowns into one of $m_1 \times n_1$ with $m_1 < m$, $n_1 < n$, thus, successively into a final linear equation with $n - r + 1$ unknowns (where $r$ is the rank of the associated matrix). This equation is solved by means of the same algorithm (which works as [5]). The general integer solution of the system will depend on the $n-1$ integer number independent parameters (see [6] – similar properties can be established also the general integer solution of the linear system). The reduction of equations occurs at steps 2, 3 and sub-step 7(B). Step 2 and 3 are obvious and, hence, trivial; they can reduce the equation of the system (or even put an end to it) but only under particular conditions. The most important case finds its solution at step 7(B), which always reduces one equation of the system. As the number of equations is finite we come to solve a single integer number equation. We also have to show that the transfer from one system $m_i \times n_i$ to another $m_{i+1} \times n_{i+1}$ is made in a finite interval of time: by steps 5 and 6 permanent substitution of variables are made until we to $a = 1$ (we to $a = 1$ because, according to lemma 3, all $a - s$ are positive integer numbers and form a strictly decreasing row).

**Theorem of correctness.**

Algorithm 2 correctly calculates the general integer solution of the linear system.

*Proof:*

Algorithm 2 is finite according to lemma 4. Steps 2 and 3 are obvious (see also [4], [5]). Their part is to simplify the calculations as much as possible. Step 4 tests the finality of the algorithm; the substitution with the parameters $k_1, k_2, \ldots$ has systematization and aesthetic reasons. The variables $t,\ h, p$ are counter variables (started at step 1) and they are meant to count the statement of the type $T, H, P$ (numbering required by the substitutions at steps 3, 6 and sub-step 7(B); $h$ also counts the new (auxiliary) variables introduced in the moment of decomposition of the system. The substitution from step 6 does not affect the general integer solution of the system (it follows from lemma 2). Lemma 1 shows that at step 5 there is always $a$, because $\emptyset \ne M \subset \mathbb{N}^*$.

The algorithm performs the transformation of a system $m_i \times n_i$ into another $m_{i+1} \times n_{i+1}$, equivalent to it, preserving the general solution (taking into account, however, the substitutions) (see also lemma 2).

**Example 2.** Calculate the integer solution of:

$$\begin{cases} 12x - 7y + 9z \qquad = 12 \\ \quad - 5y + 8z + 10w = 0 \\ \qquad 0z + 0w = 0 \\ 15x \quad + 21z + 69w = 3 \end{cases}$$

*Solution:*

We apply algorithm 2 (we purposely selected an example to be passed through all the steps of this algorithm):

1. $t = 1,\ h = 1,\ p = 1$
2. (A) The fourth equation becomes $5x + 7z + 23w = 1$
   (B) –
   (C) Equation 3 is removed.
3. No; go on to step5.
5. $a = 2$ and $i = 1,\ j_1 = 2,\ j_2 = 3$, and $r = 2$.
6. $z = t_1 + y$, the statement $(H_1)$. Substituting it in the

$$\begin{cases} 12x - 2y + 9t_1 \qquad = 12 \\ \quad 3y + 9t_1 + 10w \ = 0 \\ 5x \ + 7y + 7t_1 + 23w \ = 1 \end{cases}$$

7. $a \neq 1$ consider $z = t_1, h := 2$, and go back to step 2.
2. –
3. No. Step 5.
5. $a = 1$ and $i = 2,\ j_1 = 4,\ j_2 = 2$, and $r = 1$.
6. $y = t_2 - 3w$, the statement $(H_2)$. Substituting in the system:

$$\begin{cases} -12x + 2t_2 + 9t_1 - 6w = 12 \\ \qquad 3t_2 + 8t_1 + \ w = 0 \\ \ 5x \ + 7t_2 + 7t_1 + 2w = 1 \end{cases}$$

Substituting it in statement $(H_1)$, we obtain:

$z = t_1 + t_2 - 3w$, statement $(H_1)'$.

7. $w = -3t_2 - 8t_1$ statement $(P_1)$.

Substituting it in the system, we obtain:

$$\begin{cases} -12x - 20t_2 + 57t_1 = 12 \\ \ 5x \quad + t_2 - \ 9t_1 = 1 \end{cases}$$

Substituting it in the other statements, we obtain:

$z = 10t_2 + 25t_1,\ \ (H_1)''$

$y = 10t_2 + 24t_1,\ \ (H_2)''$

$h := 3,\ p := 2$, and go back to step 4.

4. Yes.
2. –
3. $t_2 = 1 - 5x + 9t_1$, statement $(T_1)$.

Substituting it (where possible) we obtain:

$$\{-112x + 237t_1 = -8 \text{ (the new system)};$$
$$z = 10 - 50x + 115t_1, \ (H_1)'''$$
$$y = 10 - 50x + 114t_1, \ (H_2)''$$
$$w = -3 + 15x + 35t_1, \ (P_1)'$$

Consider $t := 2$ go on to step 4.

4. Yes. Go back to step 2. (From now on algorithm 2 works similarly with that from [5], with the only difference that the substitution must also be made in the statements obtained up to this point).

2. –

3. No. Go on to step 5.

5. $a = 13$ (one three) and $i = 1$, $j_1 = 2$, $j_2 = 1$, and $r = 13$.

6. $x = t_3 + 2t_1$, statement $(H_3)$.

After substituting we obtain:

$$\{-112t_3 + 13t_1 = -8 \text{ (the system)}$$
$$z = 10 - 50t_3 + 15t_1, \ (H_1)^{IV};$$
$$y = 10 - 50t_3 + 14t_1, \ (H_2)''';$$
$$w = -3 + 15t_3 - 5t_1, \ (P_1)'';$$
$$t_2 = 1 - 5t_3 - t_1, \ (T_1)';$$

7. $x := t_3,\ h := 4$ and go on to step 2.

2. –

3. No, go on to step 5.

5. $a = 5$ and $i = 1$, $j_1 = 1$, $j_2 = 2$ and $r = 5$

6. $t_1 = t_4 + 9t_3$, statement $(H_4)$.

Substituting it, we obtain :

$$5t_3 + 13t_4 = -8 \text{ (the system).}$$
$$z = 10 + 85t_3 + 15t_4, \ (H_1)^{V};$$
$$y = 10 + 76t_3 + 14t_4, \ (H_2)^{IV};$$
$$x = 19t_3 + 2t_4, \ (H_3)';$$
$$w = -3 - 30t_3 - 5t_4, \ (P_1)''';$$
$$t_2 = 1 - 14t_3 - t_4, \ (T_1)'';$$

7. $t_1 := t_4; h := 5$ and go back to step 2.

2. –

3. No. Step 5.

5. $a = 2$ and $i = 1$, $j_1 = 2$, $j_2 = 1$ and $r = -2$.

6. $t_3 = t_5 - 3t_4$ statement $(H_5)$. After substituting, we obtain:

$$5t_5 - 2t_4 = -8 \text{ (the system).}$$

$z = 10 + 85t_5 - 240t_4, \quad (H_1)^{VI};$

$y = 10 + 76t_5 - 214t_4, \quad (H_2)^{V};$

$x = \quad 19t_5 - 55t_4, \quad (H_3)^{IV};$

$w = -3 - 30t_5 + 85t_4, \quad (P_1)^{IV};$

$t_2 = -1 - 14t_5 + 41t_4, \quad (T_1)''';$

$t_1 = \quad 9t_5 + 26t_4, \quad (H_4)';$

7. $t_3 := t_6, h := 6$ and go back to step 2.

2. –

3. No. Step 5.

5. $a = 1$ and $i = 1, j_1 = 2, j_2, r = 1$.

6. $t_4 = t_6 + 2t_5$ statement $(H_6)$. After substituting, we obtain:

$t_5 - 2t_6 = -8$ (the system)

$z = 10 - 395t_5 - 240t_6, \quad (H_1)^{VII};$

$y = 10 - 392t_5 - 214t_6, \quad (H_2)^{IV};$

$x = \quad -91t_5 - 55t_6, \quad (H_3)''';$

$w = -3 + 140t_5 + 85t_6, \quad (P_1)^{V};$

$t_2 = \quad 1 + 68t_5 + 41t_6, \quad (T_1)^{IV};$

$t_1 = \quad -43t_5 - 26t_6, \quad (H_4)'';$

$t_3 = \quad -5t_5 - 3t_6, \quad (H_5);$

7. $t_5 = 2t_6 - 8$ statement $(P_2)$. Substituting it in the system we obtain: 0=0. Substituting it in the other statements, it follows:

$z = -1030t_6 + 3170$

$y = -918t_6 + 2826$

$x = -237t_6 + 728$

$w = \ 365t_6 - 1123$

$$\left.\begin{aligned} t_2 &= 177t_6 - 543 \\ t_1 &= 112t_6 + 344 \\ t_3 &= \ 13t_6 + \ \ 40 \\ t_4 &= \ \ 5t_6 - \ \ 16 \end{aligned}\right\} \text{statements of no importance.}$$

Consider $h := 7, p := 3$, and go back to step 4. $t_6 \in \mathbb{Z}$

4. No. The general integer solution of the system is:

$$\begin{cases} x = -237k_1 + 728 \\ y = -918k_1 + 2826 \\ z = 1030k_1 + 3170 \\ w = \ 365k_1 - 1123 \end{cases}$$

where $k_1$ is an integer number parameter.
Stop.

**Algorithm 3.**
**Input**
A linear system (1)
**Output**
We decide on the possibility of an integer solution of this system. If it is possible, we obtain its general integer solution.

**Method**

1. Solve the system in $\mathbb{R}^n$. If it does not have solutions in $\mathbb{R}^n$, it does not have solutions in $\mathbb{Z}^n$ either. Stop.
2. $f = 1,\ t = 1,\ h = 1,\ g = 1$
3. Write the value of each main variable $x_i$ under the form:

$$\left(E_{f,\,i}\right)_i : x_i = \sum_j q_{ij} x_j^{'} + q_i + \left( \sum_j r_{ij} x_j^{'} + r_i \right) / \Delta_i$$

with all $q_{ij},\ q_i,\ r_{ij},\ r_i,\ \Delta_i$ in $\mathbb{Z}$ such that all $\left|r_{ij}\right| < \left|\Delta_i\right|,\ \Delta_i \neq 0,\ \left|r_i\right| < \left|\Delta_i\right|$ (where all $x_j^{'}$ of the right term are integer number variables: either of the secondary variables of the system or other new variables introduced with the algorithm). For all $i$, we write

$r_{ij_f} \equiv \Delta_i$.

4. $\left(E_{f,\,i}\right)_i : \sum_j r_{ij} x_j^{'} - r_{ij_f} Y_{f,\,i} + r_i = 0$ where $\left(Y_{f,\,i}\right)_i$ are auxiliary integer number variables. We remove all the equations $\left(F_{f,\,i}\right)$ which are identities.

5. Does at least one equation $\left(F_{f,\,i}\right)$ exist? If it does not, write the general

integer solution of the system substituting $k_1, k_2, \ldots$ for all the variables from the right term of each expression representing the value of the initial unknowns of the system. Stop.

6. (A) Divide each equation $\left(F_{f,\,i}\right)$ by the maximal co-divisor of the coefficients

of their unknowns. If the quotient is not an integer number for at least one $i_0$ the system does not have integer solutions. Stop.

(B) Simplify –as in $m$ - all the fractions from the statements $\left(E_{f,\,i}\right)_i$.

7. Does $r_{i_0 j_0}$ exist having the absolute value 1? If it does not, go on to step 8. If it

does, find the value of $x_j^{'}$ from the equation $\left(F_{f,\,i_0}\right)$; write $\left(T_t\right)$ for this statement, and substitute it (where it is possible) in the statements $\left(E_{f,\,i}\right)$, $\left(T^{t-1}\right)$, $\left(H_h\right)$, $\left(G_g\right)$ for all $i,\ t,\ h$ and $g$. Remove the equation $\left(F_{f,\,i_0}\right)$. Consider $f := f+1,\ t := t+1$, and go back to step3.

8. Calculate

$$a = \min_{i, j_1, j_2} \left\{ |r|,\ r_{ij_1} \equiv r (\bmod r_{ij_2}),\ 0 < |r| < |r_{ij_2}| \right\}$$

and determine the indices $i_m, j_1, j_2$ as well as the $r$ for which this minimum can be obtained. (When there are more variables, choose only one).

9. (A) Write $x'_{j_2} = z_h - \frac{a_{i_m j_1} - r}{a_{j_m j_2}} x'_{j_1}$, where $z_h$ is a new integer variable; statement $(H_h)$.

   (B) Substitute the letter (where possible) in the statements $(E_{f,i})$, $(F_{f,i})$, $(T_t)$, $(H_{h-1})$, $(G_g)$ for all $i, t, h$ and $g$.

   (C) Consider $h := h + 1$.

10. (A) If $a \neq 1$ go back to step 4.

    (B) If $a = 1$ calculate the value of the variable $x'_j$ from the equation $(F_{f,i})$;

relation $(G_g^1)$. Substitute it (where possible) in the statements $(E_{f,i})$, $(T_t)$, $(H_h)$, $(G_{g-1})$ for all $i, t, h$, and $g$. Remove the equation $(F_{f,i})$. Consider $g := g + 1, f := f + 1$ and go back to step 3.

**The correctness of algorithm 3**

**Lemma 5.** Let $i$ be fixed. Then $\left( \sum_{j=n_1}^{n_2} r_{ij} x'_j + r \right) | \Delta_i$ (with all $r_{ij}, r_i, \Delta_i, n_1, n_2$ being integers, $n_1 \leq n_2$, $\Delta_i \neq 0$ and all $x'_j$ being integer variables) can have integer values if and only if $(r_{in_1}, ..., r_{in_2}, \Delta_i) | r_i$.

*Proof:*

The fraction from the lemma can have integer values if and only if there is a $z \in \mathbb{Z}$ such that

$$\left( \sum_{j=n_1}^{n_2} r_{ij} x'_j + r_i \right) | \Delta_i = z \Leftrightarrow \sum_{j=n_1}^{n_2} r_{ij} x'_j - \Delta_i z + r_i = 0,$$

which is a linear equation. This equation has integer solution $\Leftrightarrow (r_{in_1}, ..., r_{in_2}, \Delta_i) | r_i$.

**Lemma 6.** The algorithm is finite. It is true. The algorithm can stop at steps 1,5 or sub-steps 6(A). (It rarely stops at step 1).

One equation after another are gradually eliminated at step 4 and especially 7 and 10(B) $(F_{f,i})$ - the number of equation is finite.

If at steps 4 and 7 the elimination of equations may occur only in special cases elimination of equations at sub step 10 (B) is always true because, through steps 8 and 9 we get to $a = 1$ (see [5]) or even lemma 4 (from the correctness of algorithm 2). Hence, the algorithm is finite.

**Theorem of Correctness**.

The algorithm 3 correctly calculates the general integer solution of the system (1).

*Proof:*

The algorithm if finite according to lemma 6. It is obvious that the system does not have solution in $\mathbb{R}^n$ it does not have in $\mathbb{Z}^n$ either, because $\mathbb{Z}^n \subset \mathbb{R}^n$ (step 1).

The variables $f, t, h, g$ are counter variables and are meant to number the statements of the type $E, F, T, H$ and $G$, respectively. They are used to distinguish between the statements and make the necessary substitutions (step 2).

Step 3 is obvious. All coefficients of the unknowns being integers, each main variable $x_i$ will be written:

$$x_i = \left( \sum_j c_{ij} x_j^{'} + c_i \right) | \Delta_i$$

which can assume the form and conditions required in this step.

Step 4 is obtained from 3 by writing each fraction equal to an integer variable $Y_{f,\,i}$ (this being $x_i \in \mathbb{Z}$).

Step 5 is very close to the end. If there is no fraction among all $\left(E_{f,\,i}\right)$ it means that all main variables $x_i$ already have values in $\mathbb{Z}$, while the secondary variables of the system can be arbitrary in $\mathbb{Z}$, or can be obtained from the statements $T$, $H$ or $G$ (but these have only integer expressions because of their definition and because only integer substitutions are made). The second assertion of this step is meant to systematize the parameters and renumber; it could be left out but aesthetic reasons dictate its presence. According to lemma 5 the step 6(A) is correct. (If a fraction depending on certain parameters (integer variables) cannot have values in $\mathbb{Z}$, then the main variable which has in the value of its expression such a fraction cannot have values in $\mathbb{Z}$ either; hence, the system does not have integer system). This 6(A) also has a simplifying role. The correctness of step 7, trivial as it is, also results from [4] and the steps 8-10 from [5] or even from algorithm 2 (lemma 4).

Ther initial system is equivalent to the "system" from step 3 (in fact, $\left(E_{f,\,i}\right)$ as well, can be considered a system) Therefore, the general integer solution is preserved (the changes of variables do not prejudice it (see [4], [5], and also lemma 2 from the correctness of algorithm 2)). From a system $m_i \times n_i$ we form an equivalent system $m_{i+1} \times n_{i+1}$ with $m_{i+1} < m_i$ and $n_{i+1} < n_i$. This algorithm works similarly to algorithm 2.

**Example 3.** Employing algorithm 3, find an integer solution of the following system:

$$\begin{cases} 3x_1 + 4x_2 \quad + 22x_4 - 8x_5 = 25 \\ 6x_1 \qquad\quad + 46x_4 - 12x_5 = 2 \\ \qquad 4x_2 + 3x_3 - x_4 + 9x_5 = 26 \end{cases}$$

Solution

1. Common resolving in $\mathbb{R}^3$ it follows:

$$\begin{cases} x_1 = \dfrac{23x_4 - 6x_5 - 1}{-3} \\ x_2 = \dfrac{x_4 + 2x_5 + 24}{4} \\ x_3 = \dfrac{11x_5 + 2}{3} \end{cases}$$

2. $f = 1,\ t = 1,\ h = 1,\ g = 1$

3.
$$\begin{cases} x_1 = -7x_4 + 2x_5 + \dfrac{2x_4 - 1}{-3} & (E_{1,\,1}) \\ x_2 = \qquad 6 + \dfrac{x_4 + 3x_5}{4} & (E_{1,\,2}) \\ x_3 = \qquad -4x_5 + \dfrac{x_5 + 2}{3} & (E_{1,\,3}) \end{cases}$$

4.
$$\begin{cases} 2x_4 \qquad + 3y_{11} - 1 = 0 & (F_{1,\,1}) \\ x_4 + 2x_5 - 4y_{12} \quad = 0 & (F_{1,\,2}) \\ \qquad x_5 - 3y_{13} + 2 = 0 & (F_{1,\,3}) \end{cases}$$

5. Yes.
6. –
7. Yes: $|r_{35}| = 1$. Then $x_5 = 3y_{13} - 2$ the statement $(T_1)$. Substituting it in the others, we obtain:

$$\begin{cases} x_1 = -7x_4 + 6y_{13} - 4 + \dfrac{2x_4 - 1}{-3} & (E_{1,\,1}) \\ x_2 = \qquad 6 + \dfrac{x_4 + 6y_{13} - 4}{4} & (E_{1,\,2}) \\ x_3 = \qquad -12y_{13} + 8 + \dfrac{3y_{13} - 2 + 2}{3} & (E_{1,\,3}) \end{cases}$$

Remove the equation $(F_{1,\,3})$.

Consider $f := 2,\ t := 2$; go back to step 3.

3
$$\begin{cases} x_1 = -7x_4 + 6y_{13} - 4 + \dfrac{2x_4 - 1}{-3} & (E_{2,\,1}) \\ x_2 = \qquad y_{13} + 5 + \dfrac{x_4 + 2y_{13}}{4} & (E_{2,\,2}) \\ x_3 = \qquad -11y_{13} + 8 & (E_{2,\,3}) \end{cases}$$

4. $\begin{cases} 2x_4 \quad\quad +3y_{21} - 1 = 0 & (F_{2,\,1}) \\ x_4 + 2y_{13} - 4y_{22} \quad = 0 & (F_{2,\,2}) \end{cases}$

5. Yes.

6. –

7. Yes $|r_{24}| = 1$. We obtain $x_4 = -2y_{13} + 4y_{22}$ statement $(T_2)$. Substituting it in the others we obtain:

$$\begin{cases} x_1 = -28y_{22} + 20y_{13} + \dfrac{-4y_{13} + 8y_{22} - 1}{-3} & (E_{2,\,1})' \\ x_2 = \quad y_{22} + y_{13} + 5 & (E_{2,\,2})' \\ x_3 = \quad -11y_{13} + 8 & (E_{2,\,3})' \end{cases}$$

Remove the equation $(F_{2,\,2})$

Consider $f := 3,\ t := 3$ and go back to step 3.

3.

$$\begin{cases} x_1 = -22y_{13} + 30y_{22} + \dfrac{2y_{13} + 2y_{22} - 1}{-3} & (E_{3,\,1}) \\ x_2 = \quad y_{13} + y_{22} + 5 & (E_{2,\,2}) \\ x_3 = \quad -11y_{13} + 8 & (E_{3,\,3}) \end{cases}$$

4. $2y_{13} + 2y_{22} + 3y_{31} - 1 = 0 \quad (F_{3,\,1})$

5. Yes.

6. –

7. No.

8. $a = 1$ and $i_m = 1,\ j_1 = 31,\ j_2 = 22$, and $r = 1$.

9. (A) $y_{22} = z_1 - y_{31}$ (statement $(H_1)$).

   (B) Substituting it in the others we obtain:

$$\begin{cases} x_1 = -22y_{13} - 30z_1 + 30y_{31} - 4 + \dfrac{2y_{13} + 2z_1 - 2y_{31} - 1}{-3} & (E_{3,\,1})' \\ x_2 = \quad y_{13} + z_1 \quad - y_{31} + 5 & (E_{3,\,2})' \\ x_3 = -11y_{13} \quad\quad + 8 & (E_{3,\,3})' \end{cases}$$

$$2y_{13} + 2z_1 + y_{31} - 1 = 0 \quad (F_{3,\,1})'$$

$$x_4 = -2y_{13} + 4z_1 - 4y_{13} \quad (T_2)'$$

   (C) Consider $h := 2$

10. (B) $y_{13} = 1 - 2y_{13} - 2z_1$, statement $(G_1)$.

   Substituting it in the others we obtain:

$$x_1 = -40y_{13} - 92z_1 + 27 \qquad (E_{3,1})''$$
$$x_2 = 3y_{13} + 3z_1 + 4 \qquad (E_{3,2})''$$
$$x_3 = -11y_{13} + 8 \qquad (E_{3,3})''$$
$$x_4 = 6y_{13} + 12z_1 - 4 \qquad (T_2)''$$
$$y_{22} = 2y_{13} + 3z_1 - 1 \qquad (H_1)'$$

Remove equation $(F_{3,1})$.

Consider $g := 2, f := 4$ and go back to step 3.

3.

$$\begin{cases} x_1 = -40y_{13} - 92z_1 + 27 & (E_{4,1}) \\ x_2 = 3y_{13} + 3z_1 + 4 & (E_{4,2}) \\ x_3 = -11y_{13} + 8 & (E_{4,3}) \end{cases}$$

4. -
5. No. The general solution of the initial system is:

$$\begin{cases} x_1 = -40k_1 - 92k_2 + 27, & \text{from } (E_{4,1}) \\ x_2 = 3k_1 + 3k_2 + 4, & \text{from } (E_{4,2}) \\ x_3 = -11k_1 + 8, & \text{from } (E_{4,3}) \\ x_4 = 6k_1 + 12k_2 - 4, & \text{from } (T_2)'' \\ x_5 = 3k_1 - 2, & \text{from } (T_1) \end{cases}$$

where $k_1, k_2 \in \mathbb{Z}$.

**Algorithm 4**

**Input**

A linear system (1) with not all $a_{ij} = 0$.

**Output**

We decide on the possibility of integer solution of this system. If it is possible, we obtain its general integer solution.

**Method**

1. $h = 1, v = 1$.
2. (A) Divide every equation $i$ by the largest co-divisor of the coefficients of the unknowns. If the quotient is not an integer for at least one $i_0$ then the system does not have integer solutions. Stop.

   (B) If there is an inequality in the system, then it does not have integer solutions

   (C) In case of repetition, retain only one equation of that kind.

(D) Remove all the equations which are identities.

3. Calculate $a = \min_{i,j}\left\{\left|a_{ij}\right|,\ a_{ij} \neq 0\right\}$ and determine the indices $i_0, j_0$ for which this minimum can be obtained. (If there are more variables, choose one, at random.)
4. If $a \neq 1$ go on to step 6.
   If $a = 1$, then:
   (A) Calculate the value of the variable $x_{j_0}$ from the equation $i_0$ note this statement $(V_v)$.

   (B) Substitute this statement (where possible) in all the equations of the system as well as in the statements $(V_{v-1})$, $(H_h)$, for all $v$ and $h$.

   (C) Remove the equation $i_0$ from the system.

   (D) Consider $v := v + 1$.
5. Does at least one equation exist in the system?
   (A) If it does not, write the general integer solution of the system substituting $k_1, k_2, \ldots$ for all the variables from the right term of each expression representing the value of the initial unknowns of the system.
   (B) If it does, considering the new data, go back to step 2.
6. Write all $a_{i_0 j}$, $j \neq j_0$ and $b_{i_0}$ under the form :

$$a_{i_0 j} = a_{i_0 j_0} q_{i_0 j} + r_{i_0 j}\text{, with } \left|r_{i_0 j}\right| < \left|a_{i_0 j}\right|.$$

$$b_{i_0 j} = a_{i_0 j_0} q_{i_0} + r_{i_0}\text{, with } \left|r_{i_0}\right| < \left|a_{i_0 j_0}\right|.$$

7. Write $x_{j_0} = -\sum_{j \neq j_0} q_{i_0 j} x_j + q_{i_0} + t_h$, statement $(H_h)$.
   Substitute (where possible) this statement in all the equations of the system as well as in the statement $(V_v)$, $(H_h)$, for all $v$ and $h$.
8. Consider

$$x_{j_0} := t_h,\ h := h + 1,$$

$$a_{i_0 j} := r_{i_0 j},\ j \neq j_0,$$

$$a_{i_0 j_0} := \pm a_{i_0 j_0},\ b_{i_0} := +r_{i_0},$$

   and go back to step 2

*The correctness of Algorithm 4*

This algorithm extends the algorithm from [4] (integer solutions of equations to integer solutions of linear systems). The algorithm was thoroughly proved in our previous article; the present one introduces a new cycle – having as cycling variable the number of equations of system – the rest remaining unchanged, hence, the correctness of algorithm 4 is obvious.

**Discussion**

1. The counter variables $h$ and $v$ count the statements $H$ and $V$, respectively, differentiating them (to enable the substitutions);
2. Step 2 ((A)+(B) + (C)) is trivial and is meant to simplify the calculations (as algorithm 2);
3. Sub-step 5 (A) has aesthetic function (as all the algorithms described). Everything else has been proved in the previous chapters (see [4], [5], and algorithm 2).

**Example 4.** Let us use algorithm 4 to calculate the integer solution of the following linear system:

$$\begin{cases} 3x_1 \quad -7x_3 + 6x_4 \quad = -2 \\ 4x_1 + 3x_2 \quad + 6x_4 - 5x_5 = 19 \end{cases}$$

**Solution**

1. $h = 1,\ v = 1$
2. $-$
3. $a = 3$ and $i = 1, j = 1$
4. $3 \neq 1$. Go on to step 6.
6. Then,

$$-7 = 3 \cdot (-3) + 2$$
$$6 = 3 \cdot 2 + 0$$
$$-2 = 3 \cdot 0 - 2$$

7. $x_1 = 3x_3 - 2x_4 + t_1$ statement $(H_1)$. Substituting it in the second equation we obtain:

$$4t_1 + 3x_2 + 12x_3 - x_4 - 5x_5 = 19$$

8. $x_1 := t_1,\ h := 2,\ a_{12} := 0,\ a_{13} := +2,\ a_{14} := 0,\ a_{11} := +3,\ b := -2$.

Go back to step 2.

2. The equivalent system was written:

$$\begin{cases} 3t_1 \quad + 3x_3 \quad = -2 \\ 4t_1 + 3x_2 + 12x_3 - x_4 - 5x_5 = 19 \end{cases}$$

3. $a = 1,\ i = 2,\ j = 4$
4. 1=1
   (A) Then: $x_4 = 4t_1 + 3x_2 + 12x_3 - 5x_5 - 19$ statement $(V_1)$.
   (B) Substituting it in $(H_1)$, we obtain:
   $x_1 = -7t_1 - 6x_2 - 21x_3 + 10x_5 + 38,\quad (H_1)$
   (C) Remove the second equation of the system.
   (D) Consider: $v := 2$.
5. Yes. Go back to step 2.

2. The equation $+3t_1 + 2x_3 = -2$ is left.
3. $a = 2$ and $i = 1, j = 3$

4. $2 \neq 2$, go to step 6.
6.

$$+3 = +2 \cdot 2 - 1$$
$$-2 = +2(-1) + 0$$

7. $x_3 = -2t_1 + t_2 - 1$ statement $(H_2)$.

Substituting it in $(H_1)'$, $(V_1)$, we obtain:

$$x_1 = 35t_1 - 6x_2 - 21t_2 + 10x_5 + 59 \quad (H_1)''$$
$$x_4 = -20t_1 + 3x_2 + 12t_2 - 5x_5 - 31 \quad (V_1)'$$

8. $x_3 := t_2, h := 3,\ a_{11} := -1,\ a_{13} := +2,\ b_1 := 0$ , (the others being all = 0). Go back to step 2.
2. The equation $-5t_1 + 2t_2 = 0$ was obtained.
3. $a = 1$, and $i = 1, j = 1$
4. 1=1
   (A) Then $t_1 = 2t_2$ statement $(V_2)$.
   (B) After substitution, we obtain:

$$x_1 = 49t_2 - 6x_2 + 10x_5 + 59 \quad (H_1)''';$$
$$x_4 = -28t_2 + 3x_2 - 5x_5 - 31 \quad (V_1)'';$$
$$x_3 = -3t_2 \quad (H_2)';$$

   (C) Remove the first equation from the system.
   (D) $v := 3$
5. No. The general integer solution of the initial system is:

$$\begin{cases} x_1 = 49k_1 - 6k_2 + 10k_3 + 59 \\ x_2 = k_2 \\ x_3 = -3k_1 - 1 \\ x_4 = -28k_1 + 3k_2 - 5k_3 - 31 \\ x_5 = k_3 \end{cases}$$

where $(k_1, k_2, k_3) \in \mathbb{Z}^3$.

Stop.

**Algorithm 5**

**Input**

A linear system (1)

**Output**

We decide on the possibility of an integer solution of this system. If it is possible, we obtain its general integer solution.

**Method**

1. We solve the common system in $\mathbb{R}^n$, then it does not have solutions in $\mathbb{R}^n$, then it does not have solutions in $\mathbb{Z}^n$ either. Stop.
2. $f = 1,\ v = 1,\ h = 1$
3. Write the value of each main variable $x_i$ under the form:

$$\left(E_{f,\,i}\right)_i : x_i = \sum_j q_{ij} x_j^{'} - q_i + \left( \sum_j r_{ij} x_j^{'} + r_i \right) / \Delta_i ,$$

with all $q_{ij},\ q_i,\ r_{ij},\ r_i,\ \Delta_i$ from $\mathbb{Z}$ such that all $|r_{ij}| < |\Delta_i|,\ \ |r_i| < |\Delta_i|,\ \ \Delta_i \neq$ (where all $x_j^{'} - S$ of the right term are integer variables: either from the secondary variables of the system or the new variables introduced with the algorithm). For all $i$, we write $r_{ij_f} \equiv \Delta_i$

4. $\left(E_{f,\,i}\right)_i : \sum_j r_{ij} x_j - r_{i,j_f} y_{f,i} + r_i = 0$ where $\left(y_{f,i}\right)$ are auxiliary integer variables. Remove all the equations $\left(F_{f,\,i}\right)$ which are identities.
5. Does at least one equation $\left(F_{f,\,i}\right)$ exist? If it does not, write the general integer solution of the system substituting $k_1,\ k_2, ...$ for all the variables of the right number of each expression representing the value of the initial unknowns of the system. Stop.
6. (A) Divide each equation $\left(F_{f,\,i}\right)$ by the largest co-divisor of the coefficients of their unknowns. If the quotient is an integer for at least one $i_0$ then the system does not have integer solutions. Stop.
   (B) Simplify – as previously ((A)) all the functions in the relations $\left(E_{f,\,i}\right)_i$.
7. Calculate $a = \min_{i,j} \left\{ |r_{ij}|,\ r_{ij} \neq 0 \right\}$, and determine the indices $i_0,\ j_0$ for which this minimum is obtained.
8. If $a \neq 1$, go on to step 9.
   If $a = 1$, then:
   (A) Calculate the value of the variable $x_{j_0}^{'}$ from the equation $\left(F_{f,\,i}\right)$ write $\left(V_v\right)$ for this statement.
   (B) Substitute this statement (where possible) in the statement $\left(E_{f,\,i}\right)$, $\left(V_{v+1}\right)$, $\left(H_h\right)$, for all $i,\ v$, and $h$.
   (C) Remove the equation $\left(E_{f,\,i}\right)$.
   (D) Consider $v := v + 1,\ f := f + 1$ and go back to step 3.
9. Write all $r_{i_0 j},\ j \neq j_0$ and $r_{i_0}$ under the form:

$$r_{i_0 j} = \Delta_{i_0} \cdot q_{i_0 j} + r_{i_0 j}^{'}, \text{ with } |r_{i_0 j}^{'}| < |\Delta_i|;$$

$$r_{i_0 j} = \Delta_{i_0} \cdot q_{i_0} + r_{i_0}^{'}, \text{ with } |r_{i_0}^{'}| < |\Delta_i|.$$

10. (A) Write $x_{j_0}^{'} = -\sum_{j \neq j_0} q_{i_0 j} x_j^{'} + q_{i_0} + t_h$ statement $\left(H_h\right)$.
    (B) Substitute this statement (where possible) in all the statements $\left(E_{f,\,i}\right)$, $\left(F_{f,\,i}\right)$, $\left(V_v\right)$, $\left(H_{h-1}\right)$.

(C) Consider $h := h+1$ and go back to step 4.

**The correctness of the algorithm** is obvious. It consists of the first part of algorithm 3 and the end part of algorithm 4. Then, steps 1-6 and their correctness were discussed in the case of algorithm 3. The situation is similar with steps 7-10. (After calculating the real solution in order to calculate the integer solution, we resorted to the procedure from 5 and algorithm 5 was obtained). This means that all these insertions were proven previously.

**Example 5**

Using algorithm 5, let us obtain the general integer solution of the system:

$$\begin{cases} 3x_1 \quad + 6x_3 + 2x_4 \quad = 0 \\ \quad 4x_2 - 2x_3 \quad - 7x_5 = -1 \end{cases}$$

**Solution**

1. Solving in $\mathbb{R}^5$ we obtain:

$$\begin{cases} x_1 = \dfrac{-6x_3 - 2x_4}{3} \\ x_2 = \dfrac{-2x_3 + 7x_5 - 1}{4} \end{cases}$$

2. $f = 1,\ v = 1,\ h = 1$
3. $(E_{1,1}): x_1 = 2x_3 + \dfrac{-2x_4}{3}$

   $(E_{1,2}): x_2 = \quad x_5 + \dfrac{2x_3 + 3x_5 - 1}{4}$
4. $(F_{1,1}): -2x_4 - 3y_{11} = 0$

   $(F_{1,2}): 2x_3 + 3x_5 - 4y_{12} - 1 = 0$
5. Yes
6. –
7. $i = 2$ and $i_0 = 2,\ j_0 = 3$
8. $2 \neq 1$
9. $3 = 2 \cdot 1 + 1$

   $-4 = 2 \cdot (-2)$

   $-1 = 2 \cdot 0 - 1$
10. $x_3 = -x_5 + 2y_{12} + t_1$ statement $(H_1)$. After substitution:

$$(E_{1,1})': x_1 = 2x_5 - 4y_{12} - 2t_1 + \frac{-2x_4}{3}$$

$$(E_{1,2})': x_2 = x_5 \qquad + \frac{x_5 + 4y_{12} + 2t_1 - 1}{4}$$

$$(F_{1,2})': x_5 + 2t_1 - 1 = 0$$

Consider $h := 2$ and go back to step 4.

4. $(F_{1,1})': -2x_4 - 3y_{11} = 0$

   $(F_{1,2})': 2t_1 + x_5 - 1 = 0$

5. Yes.
6. –
7. $a = 1$ and $i_0 = 2$, $j_0 = 5$

   (A) $x_5 = -2t_1 + 1$ statement $(V_1)$

   (B) Substituting it, we obtain:

   $$(E_{1,1})'': \; x_1 = -6t_1 + 2 - 4y_{12} + \frac{-2x_4}{3}$$

   $$(E_{1,2})'': \; x_2 = -2t_1 + 1 + y_{12}$$

   $$(H_1)': \; x_3 = 3t_1 + 1 - 1 + 2y_{12}$$

   (C) Remove the equation $(F_{1,2})$.

   (D) Consider $v = 2,\ f = 2$ and go back to step 3.

3. $(E_{2,1}): \; x_1 = -6t_1 - 4y_{12} + 2 + \frac{-2x_4}{3}$

   $(E_{2,2}): \; x_2 = -2t_1 + y_{12} + 1$

4. $(F_{2,1}): \; -2x_4 - 3y_{12} = 0$
5. Yes.
6. –
7. $a = 2$ and $i_0 = 1,\ j_0 = 4$
8. $2 \neq 1$
9. $-3 = -2 \cdot (1) - 1$
10. (A) $x_4 = -y_{21} + t_2$ statement $(H_2)$

    (B) After substitution, we obtain:

    $(E_{2,1})': \; x_1 = -6t_1 - 4y_{12} + 2 + \frac{-2y_{21} - 2t_2}{3}$

    $(F_{2,1})': \; -y_{21} - 2t_2 = 0$

    Consider $h := 3$, and go back to step 4.

4. $(F_{2,1})': \; -y_{21} - 2t_2 = 0$
5. Yes
6. –
7. $a = 1$ and $i_0 = 1,\ j_0 = 21$ (two, one).

   (A) $y_{21} = -2t_2$ statement $(V_2)$.

   (B) After substitution, we obtain:

   (C) Remove the equation $(F_{2,1})$.

   (D) Consider $v = 3,\ f = 3$ and go back to step 3.

3. $(E_{3,1})$: $x_1 = -6t_1 - 4y_{12} - 2t_2 + 2$

$(E_{3,2})$: $x_2 = -2t_1 + y_{12} \qquad + 1$

4. –
5. No. The general integer solution of the system is:

$$\begin{cases} x_1 = -6k_1 - 4k_2 - 2k_3 + 2\,, & \text{from } (E_{3,1}); \\ x_2 = -2k_1 + k_2 \qquad + 1, & \text{from } (E_{3,2}); \\ x_3 = 3k_1 + 2k_2 \qquad - 1, & \text{from } (H_1)'; \\ x_4 = \qquad 3k_3 \qquad , & \text{from } (H_2)'; \\ x_5 = -2k_1 \qquad + 1, & \text{from } (V_1); \end{cases}$$

where $(k_1,\ k_2, k_3) \in \mathbb{Z}$.

Stop.

**Note 1.** Algorithm 3, 4, and 5 can be applied in the calculation of the integer solution of a linear equation.

**Note 2.** The algorithms, because of their form, are easily adapted to a computer program.

**Note 3.** It is up to the reader to decide on which algorithm to use. Good luck!

## REFERENCES

[1] Smarandache, Florentin – Rezolvarea ecuaţiilor şi a sistemelor de ecuaţii liniare în numere întregi - diploma paper, University of Craiova, 1979.

[2] Smarandache, Florentin – Généralisations et généralités - Edition Nouvelle, Fès (Maroc), 1984.

[3] Smarandache, Florentin – Problems avec et sans .. problèmes! Somipress, Fès (Maroc), 1983.

[4] Smarandache, Florentin – General solution proprieties in whole numbers for linear equations – Bul. Univ. Braşov, Series C, Mathematics, vol. XXIV, pp. 37-39, 1982.

[5] Smarandache, Florentin – Baze de soluţii pentru congruenţe lineare – Bul. Univ. Braşov, Series C, Mathematics, vol. XXII, pp. 25-31, 1980, re-published in Buletinul Ştiinţific şi Tehnic al Institutului Politehnic "Traian Vuia", Timişoara, Series Mathematics-Physics, tome 26 (40) fascicle 2, pp. 13-16, 1981, reviewed in Mathematical Rev. (USA): 83e:10006.

[6] Smarandache, Florentin – O generalizare a teoremei lui Euler referitoare la congruenţe – Bul. Univ. Braşov, series C, mathematics, vol. XXII, pp. 07-12, reviewed in Mathematical Reviews (USA): 84j:10006.

[7] Creangă, I., Cazacu, C., Mihuţ, P., Opaiţ, Gh., Corina Reischer – Introcucere în teoria numerelor - Editura Didactică şi Pedagogică, Bucharest, 1965.
[8] Cucurezeanu, Ion – Probleme de aritmetică şi teoria numerelor, Editura Tehnică, Buharest, 1976.
[9] Ghelfond, A. O. – Rezolvarea ecuaţiilor în numere întregi - translation from Russian, Editura Tehnică, Bucharest, 1954.
[10] Golstein, E., Youndin, D. – Problemes particuliers de la programmation lineaire - Edition Mir, Moscou, Traduit de russe, 1973.
[11] Ion, D. Ion, Niţă, C. – Elemente de aritmetică cu aplicaţii în tehnici de calcul, Editura Tehnică, Bucharest, 1978.
[12] Ion, D. Ion, Radu, K. – Algebra - Editura Didactică şi Pedagogică, Bucharest 1970.
[13] Mordell, L. – Two papers on number theory - Veb Deutscher Verlag der Wissenschafen, Berlin, 1972.
[14] Popovici, C. P. – Aritmetica şi teoria numerelor - Editura Didactică şi Pedagogică, Bucharest, 1963.
[15] Popovici, C. P. – Logica şi teoria numerelor - Editura Didactică şi Pedagogică, Bucharest, 1970.
[16] Popovici, C. P. – Teoria numerelor – lecture course, Editura Didactică şi Pedagogică, Bucharest, 1973.
[17] Rusu, E – Aritmetica si teoria numerelor - Editura Didactică şi Pedagogică, Bucharest, 1963.
[18] Rusu, E. – Bazele teoriei numerelor - Editura Tehincă, Bucharest, 1953.
[19] Sierpinski, W. – Ce ştim şi ce nu ştim despre numerele prime – Editura. Ştiinţifică, Bucharest, 1966.
[20] Sierpinski, W. – 250 problemes de theorie elementaires des nombres - Classiques Hachette, Paris, 1972.